
\input amstex
\input amsppt.sty
\magnification=1200
\nologo

\def\proof{\demo{Proof}}
\def\endproof{\qed\enddemo}
\def\nextpart#1{\medskip{\bf (#1).}}

\def\A{{\Cal A}}
\def\C{{\check C}}
\def\D{{\bold D}}
\def\E{{\Cal E}}
\def\F{{\Cal F}}
\def\G{{\Cal G}}
\def\H{{\Cal H}}
\def\I{{\Cal I}}
\def\J{{\Cal J}}
\def\K{{\Cal K^\bullet}}
\def\KK{{\Cal K_{\<\<\sssize\infty}^\bullet}\<\<}
\def\L{{\bold L}}
\def\O{{\Cal O}}
\def\OX{{\O_{\<\<X}}}
\def\P{{\Cal P}}
\def\fP{{\text{\eighteufm P}}}
\def\Q{{\Cal Q}}
\def\R{{\bold R}}
\def\Rq#1{{\R\>Q_{#1}}}

\def\bt{{\bold t}}
\def\GG#1{{\Gamma_{\mkern-5.5mu#1}^{\phantom{.}}}}
\def\LL#1{{\L\Lambda_{#1}}}
\def\vg#1{{\varGamma_{\!\<#1}^{\phantom{\prime}}@,@,}}
\def\vgp#1{{\varGamma_{\!\<#1}'}@,@,}
\def\vgpt{{\varGamma_{\!\bt}'}\<}
\def\a{{\alpha}}

\def\Ext{{\operatorname{Ext}}}
\def\Hom{{\operatorname{Hom}}}
\def\Homb{{\operatorname{Hom}^\bullet}}
\def\sHom{{\Cal H\eurm o\eurm m}\mkern2mu}
\def\sHomb{\sHom\<\<^\bullet}
\define\RsH#1{\R\sHom\!_{#1}^\bullet}
\define\RssH#1{\R\Cal H\text{\eighteurm om}_{#1}^\bullet}

\def\Dc{{\D_{\text{\rm c}}}}

\def\qcu{^{\text{\rm qc}}}
\def\qcd{_{\text{\rm qc}}}
\def\upl{{\mkern-1.3mu\raise.2ex\hbox{$\sssize +$}}}
\def\umi{{\mkern-1.3mu\raise.2ex\hbox{$\sssize -$}}}

\def\set{\!:=}
\def\lra{\longrightarrow}
\def\iso{{\mkern8mu\longrightarrow \mkern-25.5mu{}^\sim\mkern17mu}}
\def\osi{{\mkern8mu\longleftarrow \mkern-24.5mu{}^\sim\mkern16mu}}
\def\Iso{\vbox to 0pt{\vss\hbox{$\widetilde{\phantom{nn}}$}\vskip-7pt}}
\def\drlm{\varinjlim\,}
\def\inlm{\varprojlim\,}
\def\Otimes{\underset
  {\vbox to 0pt {\vskip-1.05ex\hbox{$\sssize=$}\vss}}
     \to\otimes}
\def\smcirc{{\hbox to.7em{$\hss \eightpoint\circ \hss$}}}  
\def\circle#1{\hbox{$\bigcirc\hbox{\kern-7pt\hbox
     {{\eightpoint #1}}}\ $}}
\def\w#1{{\widetilde{#1}}}
\def\smcong{\hbox{\eightpoint\vbox to0pt
     {\vss\hbox{$\;\sim$}\vskip-2.1ex\hbox{${}={}$}}}} 
\define\upll{{\text {\smrm\char'53}}}          

\def\undercircle#1{\vbox to 0pt{\vskip.1in\circle{#1}\vss}}

\def\jog{\mkern1.5mu}
\def\>{{\mkern 1mu}}
\def\<{{\mkern-1mu}}
\def\nmb{\nomathbreak}
\def\quot#1{\mkern1.2mu\:\!#1}

%
%
\define\UnderElement#1#2#3#4{\vbox to 0pt{
\hbox{$
\llap{$\scriptstyle#1$}
\left#2\vbox to #3{}\right.
\rlap{$\scriptstyle#4$}
     $}
\vss}}

\loadeufm
\loadeurm
\loadeusm
\loadbold
\font\eighteufm=eufm8
\font\eighteurm=eurm8
\font\eightBbb=msbm8
\font\smrm=cmr5 scaled 720
\topmatter

\title {Local Homology and Cohomology on Schemes}
\endtitle

\author Leovigildo Alonso Tarr\'io, Ana Jerem\'ias L\'opez,
and Joseph Lipman \endauthor

\leftheadtext{L\.  Alonso, A\. Jerem\'{\i}as,
J\. Lipman}

\rightheadtext{Local Homology and Cohomology}

\thanks First two authors partially supported by a
Xunta de Galicia (D.O.G\. 19/11/92) travel grant. Third author partially
supported by the National Security Agency.  \endthanks

\address{Universidade de Santiago de Compostela,
Facultade de Matem\'aticas,\newline
E-15771  Santiago de Compostela, SPAIN}\endaddress

\email{lalonso\@zmat.usc.es}\endemail

\address{Universidade de Santiago de Compostela,
Facultade de Matem\'aticas,\newline
E-15771  Santiago de Compostela, SPAIN}\endaddress

\email{jeremias\@zmat.usc.es}\endemail

\address{Dept\. of Mathematics, Purdue University,
 W. Lafayette IN 47907, USA} \endaddress

\email{lipman\@math.purdue.edu}\endemail

\comment
\abstract{We prove a sheaf-theoretic derived-category generalization of
Greenlees-May duality (a far-reaching
generalization of Grothendieck's local duality theorem.)
Combined with global Grothendieck duality, this leads to a duality theorem
for proper maps of formal schemes. These results unify a scattered host of
duality theorems.}
\endabstract
\endcomment

\abstract{We prove a sheaf-theoretic derived-category generalization of
Greenlees-May duality (a far-reaching
generalization of Grothendieck's local duality theorem):  for a quasi-compact
separated scheme~$X$ and a
``proregular\kern.5pt " subscheme~$Z$---for example, any
separated noetherian scheme and any closed subscheme---there is a
sort of sheafified adjointness between local cohomology supported in ~$Z$
and left-derived completion along~$Z$.
In particular, left-derived
completion can be identified with {\it local~homology,} i.e., the homology
of~$\R\sHomb(\R\vg Z\OX,-)$.

Sheafified generalizations of a number of duality theorems scattered about
the literature result: the {\it Peskine-Szpiro duality sequence\/}
(generalizing local duality),  the {\it Warwick Duality\/} theorem of
Greenlees,  the
{\it Affine Duality\/} theorem of Hartshorne.
Using Grothendieck Duality, we also get a generalization
of a {\it Formal Duality\/} theorem of Hartshorne, and of a related
local-global duality theorem.

In a sequel we will develop the latter results further,
to study Grothendieck duality and residues on formal schemes.}
\endabstract

\endtopmatter

\document

\subheading{Introduction}
Our main result is the Duality Theorem (0.3) on a quasi-compact separated
scheme~$X$ around a {\it proregularly embedded\/} closed subscheme~$Z$.
This asserts a
sort of sheafified adjointness between
local cohomology supported in ~$Z$
and~left-derived functors of completion along~$Z$. (For complexes with
quasi-coherent homology, the precise derived-category adjoint of
local cohomology is described in Remark (0.4)(a).)
A special case---and also a
basic point in the proof---is that $(*)$:~these left-derived
completion functors can be identified with {\it local~homology,} i.e.,~the
homology of~$\R\sHomb(\R\vg Z\OX,-)$.

The technical condition ``$Z$ proregularly embedded," treated at length
in \S3, is just what is needed to make cohomology supported in ~$Z$
enjoy some good properties which are standard when $X$ is noetherian.
Indeed, it might be said that these properties hold in
the noetherian context because (as follows immediately from
the definition) {\it every closed subscheme of a
noetherian scheme is proregularly embedded.}

The assertion $(*)$ is a sheafified derived-category version of
Theorem\,2.5 in \cite{GM}.  (The particular case where $Z$ is regularly
embedded
in~$X$ had been studied, over commutative rings,
by Strebel \cite{St, pp.\;94--95, 5.9} and, in great detail,
by Matlis \cite{M2, p.\,89, Thm.\,20}.
Also, a special case of
Theorem~(0.3) appeared in \cite{Me, p.\,96} at
the beginning of the proof of 2.2.1.3.) More specifically, our
Proposition~(4.1) provides another approach to the Greenlees-May duality
isomorphism---call it $\Psi$---from local homology to left-derived completion
functors.  But this $\Psi$ is local and depends on choices, so for
globalizing there remains the non-trivial question of canonicity.  This
is dealt with in Proposition (4.2), which states that a certain natural
global map~$\Phi$ from left-derived completion functors to local
homology restricts locally to an inverse of~$\Psi$. The
map~$\Phi$ is easy to define (\S2), but we don't know any other~way to show
that it is an isomorphism.

We will exhibit in
\S5 how Theorem~(0.3) provides a unifying principle for a substantial
collection of other duality results from the literature
(listed in the introductions to those sections). For example, as
noted by Greenlees and May \cite{GM, p.\,450, Prop.\,3.8},
their theorem contains the standard Local Duality theorem of Grothendieck.
(See Remark~(0.4)(c) below for more in this vein).

\comment
Section 6 is devoted mostly to our second main result, Theorem (6.2), a form of
Grothendieck\- Duality for pseudo\kern.5pt-proper maps of noetherian formal
schemes.
Here, Theorem~(0.3) opens up the passage from ordinary to formal schemes.
\endcomment

\medskip

To describe things more precisely, we need some notation. Let $X$ be
a quasi-compact separated scheme,
let $\A(X)$ be the category of all $\OX$-modules, and let
$\A\qcd(X)\subset\A(X)$ be the full (abelian) subcategory of
quasi-coherent $\OX$-modules. The derived category
$\D(X)$ of $\A(X)$ contains full subcategories
$\D\qcd(X)\supset\nmb \D_{\text{c}}(X)$ whose objects are
the $\OX$-complexes with  quasi-coherent,
respectively coherent, homology sheaves.

Let $Z\subset X$ be a closed subset. If $X\setminus Z$ is quasi-compact then
by induction\- on~
$\min\{\>\>n \mid X\text{ can be covered by $n$ affine open subsets}\>\>\}$,
and \cite{GrD, p.\,318, (6.9.7)}, one shows that $Z$ is the
support~$\text{Supp\kern.5pt}(\OX/\I\>)$ for some
{\it finite-type\/} quasi-coherent \hbox{$\OX$-ideal~$\I$} (and conversely).
We assume throughout that $Z$ satisfies this condition.

The left-exact functor $\vg Z\:\A(X)\to\A(X)$
associates to each $\OX$-module $\F\>$ its subsheaf of sections with
support in~$Z$. We define the subfunctor $\vgp Z\subset\vg Z$ by
$$
\vgp Z\F\set\underset{{}^{n>0}}\to\drlm \,
 \sHom_{\OX}(\OX/\I^{\>n}\<,\>\jog \F\>)
  \qquad\bigl(\F\in\A(X)\bigr), \tag 0.1
$$
which depends only on ~$Z$ (not~$\I$). If $\F$ is quasi-coherent, then $\vgp
Z(\F\>)=\nmb\vg Z(\F\>)$.
The~functor $\vg Z$ (resp.~$\vgp Z$) has a right-derived functor
$\R\vg Z:\nmb\D(X)\to\D(X)$
(resp\. $\R\vgp Z\:\nmb\D(X)\to\D(X)$),
as does any functor from $\A(X)$ to an abelian category, via K-injective
resolutions \cite{Sp, p.\,138, Thm.\,4.5}.
\footnote
{See also [{\it ibid.,} p.\,133, Prop.\,3.11] or \cite{BN, \S2} for the
existence of such resolutions in module categories. (Actually, as recently
observed by Weibel, Cartan-Eilenberg
resolutions, totalized via products, will do in this case.) Moreover,
Neeman has a strikingly simple proof that hence such resolutions exist in
any abelian quotient category of a module category, i.e., by a theorem of
Gabriel-Popescu, in any abelian category---for instance $\A(X)$---with
a generator and with exact filtered $\drlm$. (Private communication.)%
}
By the universal property of derived functors, there is a
unique functorial map
$$
\gamma\:\R\vgp Z\E\to \E
$$
whose composition with $\vgp Z\E\to\R\vgp Z\E\/$ is
the inclusion map $\vgp Z\E\hookrightarrow\E$.

For proregularly embedded $Z\subset X$, {\it the derived-category map\/
$\R\vgp Z\E\to\nmb\R\vg Z\E$ induced by the inclusion\/~
$\vgp Z\hookrightarrow\nmb\vg Z$
is an isomorphism for any complex\/} $\E\in\nmb\D\qcd(X)$
(Corollary (3.2.4)). This iso\-morphism underlies the
well-known homology isomorphisms (of sheaves)
\footnote
{See \cite{H, p.\,273}, where, however, the proof seems incomplete---``way-out"
needs to begin with \cite{Gr,~p.\,22, Thm.\,6}. Alternatively, one
could use quasi-coherent injective resolutions\dots}
$$
\underset{{}^{n>0}}\to\drlm
\E{\eurm x \eurm t}^i(\OX/\I^{\>n}\<,\jog\F\>)\iso H_{\<\<Z}^i(\F\>)
\qquad (i\ge0,\ \F\in\A\qcd(X)).\tag 0.1.1
$$
\smallskip

We also consider the {\it completion functor\/} $\Lambda_Z:\A\qcd(X)\to\A(X)$
given by
$$
\nopagebreak
\Lambda_Z\F\set
 \underset{{}^{n>0}}\to\inlm\,\bigl((\OX/\I^{\>n})\otimes\F\jog\bigr)
 \qquad\bigl(\F\in\A\qcd(X)\bigr).\tag 0.2
$$
This depends only on~$Z$.
We will show in \S1 that $\Lambda_Z$ has a left-derived functor
$\L\Lambda_Z\:\D\qcd(X)\to\nmb\D(X)$, describable via flat quasi-coherent
resolutions. By the universal property of derived functors, there is a
unique functorial map
$$
\lambda\:\F\to\L\Lambda_Z\F
$$
whose composition with $\L\Lambda_Z\F\to \Lambda_Z\F\/$ is
the completion map $\F\to\nmb \Lambda_Z\F$.

\proclaim{Theorem (0.3)}
For any quasi-compact separated scheme\/~$X$
and any proregularly embedded closed subscheme\/~$Z$
$($Definition~$(3.0.1)),$\  there is a functorial isomorphism
$$
\spreadlines{1\jot}
\multline
\kern2cm\RsH{}(\R\vgp Z\E,\> \F\>)\iso
 \RsH{}(\E,\>\L\Lambda_Z \F\>)\\
  \bigl(\E\in\D(X),\ \F\in\D\qcd(X)\bigr)
\endmultline
$$
whose composition with the map\/
$\RsH{}\bigl(\E,\> \F\>\bigr)\to
 \RsH{}\bigl(\R\vgp Z\E,\> \F\>\bigr)$ induced by\/~$\gamma$ is the map\/
$\RsH{}\bigl(\E,\> \F\>\bigr)\to \RsH{}(\E,\>\L\Lambda_Z \F\>)$
induced by\/~$\lambda$.
\endproclaim

The {\it proof\/} occupies \S\S1--4; an outline is given in~\S2.
Miscellaneous corollaries and applications appear in \S5.
\medskip

{}From Theorem~(0.3) we get a commutative diagram
\newcount\SaveCatCode
\SaveCatCode=\the\catcode`@
\catcode`@=11
\expandafter\ifx\csname graphb\endcsname\relax
\alloc@4\box\chardef\insc@unt\graphb\fi
\catcode`@=\number\SaveCatCode
\setbox\graphb=\vtop{%
  \baselineskip=0pt \lineskip=0pt \lineskiplimit=0pt
  \vbox to0pt{\hbox{%
    \special{pn 8}%
    \special{pn 2}%
    \special{pa 0 0}%
    \special{pa 379 249}%
    \special{fp}%
    \special{pa 346 260}%
    \special{pa 379 249}%
    \special{fp}%
    \special{pa 376 214}%
    \special{pa 379 249}%
    \special{fp}%
    \kern  0.380in
  }\vss}%
  \kern  0.261in
}

$$
\CD
\RsH{}(\R\vgp Z\E,\>\F\>)
 @>\Iso>\vbox to0pt{\vskip-1.2ex\hbox{$\ssize\alpha$}\vss}>
  \RsH{}(\E,\>\LL Z\F\>)  \\
@VV\kern54pt{\box\graphb}V
 @V\kern-78pt\lower3ex\hbox{$\ssize\lambda'$} V\gamma' V \\ \vspace{1.5pt}
\RsH{}(\R\vgp Z\>\R\vgp Z\E,\>\F\>) @>\Iso>\phantom{.}> \RsH{}(\R\vgp Z\E,\>\LL
Z\F\>)
\endCD
$$
with horizontal isomorphisms as in~(0.3),  $\lambda'$ induced by~$\lambda$,
and $\gamma'$ induced by~$\gamma$. It~follows readily from
Lemma~(3.1.1)(2) that
the natural map $\R\vgp Z\>\R\vgp Z\E\to\nmb\R\vgp Z\E$ is an isomorphism;
hence both $\lambda'$ and~$\gamma'$ are isomorphisms, and $\alpha$ has the
{\it explicit description\/} $\alpha=\gamma'{}^{-1}\smcirc \lambda'$.
Conversely, if we knew beforehand that $\lambda'$ and
$\gamma'$ are isomorphisms, then we could {\it define\/}
$\alpha\set\gamma'{}^{-1}\smcirc \lambda'$ and recover Theorem~(0.3).
Thus we can restate the Theorem as:

\proclaim{Theorem \kern-1pt(0.3)\kern-.3pt(bis)}
\kern-2ptFor any quasi-compact separated scheme\/~$X$
and proregularly embedded closed subscheme\/~$Z,$\
the maps $\lambda$ and\/ $\gamma$ induce functorial isomorphisms
$$
\spreadlines{1.2\jot}
\multlinegap{0pt}
\multline
\qquad\RsH{}(\R\vgp Z\E,\> \F\>)\underset{\lambda'}\to\iso
 \RsH{}(\R\vgp Z\E,\> \LL Z\F\>)\underset{\gamma'}\to\osi
 \RsH{}(\E,\>\L\Lambda_Z \F\>)\\
  \bigl(\E\in\D(X),\ \F\in\D\qcd(X)\bigr).
\endmultline
$$
\endproclaim

As explained in Remark~(5.1.2), that
$\lambda'$ is an isomorphism amounts to the following Corollary.
Recall that {\it proregularity\/}
of a finite sequence $\bt\set(t_1,t_2,\dots,t_\mu)$ in a commutative
ring $A$ is defined in~(3.0.1) (where $X$ can be taken to be
$\text{Spec}(\<A))$;
and that every sequence in a noetherian ring is proregular.

\proclaim{Corollary (0.3.1)}
Let\/ $\bt$ be a proregular sequence in a commutative ring\/~$A,$\
and let $F$ be a flat\/ $A$-module, with\/ $\bt$-adic completion\/~$\widehat
F$.
Then the natural local homology maps
$H_{\bt A}^n(F)\to H_{\bt A}^n(\widehat F\>)\ \,(n\ge 0)$
are all isomorphisms.

In other words, the natural Koszul-complex map
$\KK(\bt)\otimes F\to\nmb\KK(\bt)\otimes \widehat F$ is a quasi-isomorphism
$($see $(3.1.1)(2)).$
\endproclaim

Suppose now that $X$ is affine, say $X=\text{Spec}(\<A)$, let
$\bt\set\nmb(t_1,t_2,\dots,t_\mu)$ be a proregular sequence in~$A$,
and set $Z\set\text{Spec}(\<A/\bt A)$.
With $\bt^{\<n}\set(t_1^n,\dots,t_\mu^n)$, consider the $A$-module functors
$$
\spreadlines{1\jot}
\align
\Gamma_{\!\bt}(G)&\set\underset{{}^{n>0}}\to\drlm \,
 \Hom_{A}(\<A/\bt^{\<n}\!A,\jog G), \\
\Lambda_\bt(G)&\set
 \underset{{}^{n>0}}\to\inlm\,\bigl((\<A/\bt^{\<n}\!A)\otimes G\bigr)
 \qquad(\text{$\bt$-adic completion)}.
\endalign
$$
This is the situation in \cite{GM}, and when the sequence~
$\bt$ is {\it $A$-regular,} in \cite{M2}.
The arguments used here to prove Theorem~(0.3) apply as well in the simpler
ring-theoretic context, yielding an isomorphism in the derived $A$-module
category $\D(\<A)$:\looseness=-1
$$
\quad\R\Hom_{\mkern-1.5mu A}^\bullet(\R\Gamma_{\!\bt}E,\mkern1mu F\>)\iso
 \R\Hom_{\mkern-1.5mu A}^\bullet(E,\mkern1mu\L\Lambda_\bt F\>)\qquad
  \bigl(E,F\in\D(\<A)\bigr).
\tag "$(0.3)_{\text{aff}}$"
$$

In fact $(0.3)_{\text{aff}}\>$ (with the isomorphism explicated as in (0.3)
or (0.3)bis) is essentially equivalent to~$(0.3)$ for $\E\in\D\qcd\>$,
see Remark~(d).

Suppose, for example, that $\bt$ is $A$-regular, so that there is an
isomorphism
$$
\R\Gamma_{\!\bt}(\<A)[\mu]\iso H^\mu_{\bt A}(\<A)=:K.
$$
Then for any $A$-complex~$F$, there
is a natural isomorphism \smash{$K[-\mu]\Otimes F\iso \R\Gamma_{\!\bt}(F\>)$}
(cf.~Corollary (3.2.5)), and so we have a composed isomorphism
$$
\spreadlines{1\jot}
\align
H^0\L\Lambda_\bt (F\>)
&\iso H^0\R\Hom_{\mkern-1.5mu A}^\bullet(\R\Gamma_{\!\bt}A,\>F\>)\\
&\iso H^0\R\Hom_{\mkern-1.5mu A}^\bullet(\R\Gamma_{\!\bt}A,\>
 \R\Gamma_{\!\bt}F\>)
  \iso \smash{\Hom_{\D(\<A)}\<(K,\>K\Otimes F\>)}
\endalign
$$
corresponding to the First Representation Theorem of \cite{M2, p.\,91}.
\footnote{
Matlis states the theorem for $A$-modules~$F$ which are  ``K-torsion-free"
i.e\. ([{\it ibid,} p.\,86]), the canonical map
\smash{$K\Otimes F\to K\otimes F$}\vadjust{\kern1.66pt} is an isomorphism;
and he shows for such~$F$ that the natural map
$H^0\L\Lambda_\bt (F\>)\to \Lambda_\bt (F\>)$ is an isomorphism
[{\it ibid,} p.\,89, Thm.\,21, (2)].}

\remark\nofrills{\bf Remarks (0.4).\usualspace}
{\bf(a)} Fix a quasi-compact separated scheme~$X\<$, and
write $\A$, $\A\qcd\>$, $\D$, $\D\qcd\>$, for~$\A(X)$, $\A\qcd(X)$,
$\D(X)$, $\D\qcd(X)$, respectively. Let $Z\subset X$ be a proregularly
embedded  closed subscheme. Corollary (3.2.5)(iii)
gives us the functor $\R\vg Z\:\D\qcd\to\nmb\D\qcd\>$. Theorem (0.3) yields a
{\it right adjoint\/} for this functor, as follows.\looseness =-1

The inclusion functor $\A\qcd\hookrightarrow\A$ has a right adjoint $Q$, the
``quasi-coherator''
\cite{I,~p.\,187, Lemme 3.2}. The functor $Q$, having an exact left adjoint,
preserves  K-injectivity, and it follows then  that $\R\> Q$ is right-adjoint
to the natural functor $\boldkey j\:\D(\A\qcd)\to\nmb\D$,
see \cite{Sp, p.\,129, Prop.\,1.5(b)}.
 Since $\boldkey j$~induces an {\it equivalence of categories\/}
$\D(\A\qcd)\approx\D\qcd$ (see \S1), therefore the inclusion functor
$\D\qcd\hookrightarrow\D$ has a right adjoint,
which---mildly abusing notation---we continue to
\vadjust{\goodbreak\noindent} denote by $\R\> Q$.
Thus there is a functorial isomorphism
$$
\Hom_\D\<(\E,\>\L\Lambda_Z\F)\iso
 \Hom_{\D\qcd}(\E,\> \R\> Q\>\L\Lambda_Z\F\>)
  \qquad(\E,\F\in\D\qcd).
$$
Recalling that $\R\vgp Z$ coincides with~$\R\vg Z$ on~$\D\qcd\>$,
and applying the functor $H^0\R\Gamma$ to the isomorphism in (0.3),
\footnote{
Note that $H^0\R\Gamma\>\R\H\text{\eighteurm om}^\bullet=H^0\R\Homb=\Hom_\D$,
see e.g., \cite{Sp, 5.14, 5.12, 5.17\kern.6pt}. (In order to combine left- and
right-derived functors, we must deal with {\it unbounded\/} complexes.)%
}
we deduce an {\it adjunction isomorphism}
$$
\Hom_{\D\qcd}(\R\vg Z\E,\>\F\>)\iso
 \Hom_{\D\qcd}(\E\<,\>\R\>Q\>\L\Lambda_Z\F\>)
  \qquad(\E,\F\in\D\qcd).
$$
(In this form the isomorphism doesn't sheafify, since for
open immersions $i\:U\to\nmb X$ the canonical functorial map
$i^*Q_X\to\nmb Q_U\>i^*$ is usually {\it not\/} an isomorphism.)
\smallskip
\eightpoint
For example, if $X$ is {\it affine,} say $X=\text{Spec}(\<A)$, then for any
$\Cal L\in\A(X)\>$, $Q(\Cal L)$ is the
quasi-coherent $\OX$-module $\Gamma(X,\Cal L)^\sim$
associated to the $A$-module $\Gamma(X,\Cal L)$; and hence
$$
\R\>Q(\G)\smcong
\bigl(\R\Gamma(X,\> \G)\bigr)^\sim\qquad(\G\in\D).
$$
Any complex in $\D\qcd$ is isomorphic to a K-flat quasi-coherent
$\F$ (Prop.\,(1.1)). For such an $\F$, with $\F_{/Z}$
the {\it completion\/} of~$\F$ along~$Z$ \cite{GrD, p.\,418, (10.8.2)},
 Remark (d) below, with $E=A$, implies
$$
\R\>Q\>\L\Lambda_Z\F\smcong Q\Lambda_Z\F=\bigl(\Gamma(Z,\>\F_{/Z})\bigr)^\sim.
$$
If furthermore $A$ is noetherian, $Z=\text{Spec}(\<A/I)$, and $\F\in\Dc(X)$,
then one finds, as in (0.4.1) below, that with $\hat A$ the $I$-adic completion
of~$A$,
$$
\Gamma(Z,\>\F_{/Z})\smcong \Gamma(X,\>\F\>)\otimes_{\<A}@!\hat A.
$$
\par \tenpoint
\smallskip

In more detail, Theorem~(0.3)---at least for
$\E\in\D\qcd(X)$---can be expressed via category-theoretic properties
of the endofunctors $S\set\R\vg Z$ and $T\set\R\>Q\>\LL Z$ of~$\D\qcd(X)$.
(In the commutative-ring context, use
$S\set\R\GG{\jog\bt}$ and $T\set\LL\bt\>$ instead.)

\proclaim{Theorem~$(0.3)^*$}
The canonical maps\/
$S@>\gamma>>\bold 1 @>\nu>>T\>$
$(\<\<$where $\bold 1$ is the identity functor of\/ $\D\qcd(X))$
induce functorial isomorphisms
$$
\Hom(S\E,\>S\F\>)\cong\Hom(S\E,\>\F\>)\cong
 \Hom(S\E,\>T\F\>)\cong\Hom(\E,\>T\F\>)\cong \Hom(T\E,\>T\F\>).
$$
\endproclaim

\proof (See also (5.1.1.))
The first isomorphism is given by
Lemma (0.4.2) below. The next two follow from Theorem~(0.3)(bis), giving the
adjointness of $S$ and $T\<$, as well as the isomorphism
$S\iso ST\/$ in the following Corollary. Hence:
$$
\Hom(\E,\>T\F\>)\cong\Hom(S\E,\>\F\>)
 \cong\Hom(ST\E,\>\F\>)\cong\Hom(T\E,\>T\F\>).\qquad\square
$$
\enddemo

Conversely, Theorem $(0.3)^*$, applied to arbitrary affine open subsets of~
$X\<$, yields Theorem~(0.3)(bis).

\proclaim{Corollary}
The maps $\gamma$ and $\nu$ induce functorial isomorphisms

\smallskip
\setbox1=\hbox{\rm(iii)}

$\hbox to\wd1{\rm (i)\hss}\quad S^2\iso S$.

$\hbox to\wd1{\rm (ii)\hss}\quad T\iso T^2$.

{\rm (iii)}\quad $TS\iso T$.

{\rm (iv)}\quad $S\iso ST$.
\endproclaim
\goodbreak
\proof
(i) follows, for example, from the functorial isomorphism (see (a) above)
$\Hom(S\E,\>S^2\F\>)\iso  \Hom(S\E,\>S\F\>)\>$
applied when $\E=\F$ and when $\E=S\F$.\vadjust{\kern.7pt}

(ii):  equivalent to (i) by adjunction.\vadjust{\kern1pt}

(iii): use
$\Hom(\E,\>TS\F\>)\cong\Hom(S\E,\>S\F\>)
 \cong\Hom(S\E,\>\F\>)\cong\Hom(\E,\>T\F\>).$

(iv): use
$\Hom(S\E,\>S\F\>)\cong\Hom(S\E,\>\F\>)
 \cong\Hom(S\E,\>T\F\>)\cong\Hom(S\E,\>ST\F\>).\hfill\square$
\enddemo

{\bf(b)} With notation as (a), suppose that the separated scheme~$X$ is
noetherian, so that any closed subscheme~$Z$ is proregularly embedded. On
coherent $\OX$-modules the functor~$\Lambda_Z$ is {\it exact.}
This suggests (but doesn't prove) the following concrete interpretation for the
restriction of the derived functor~$\L\Lambda_Z$ to
$\D_{\text{c}}\subset\nmb\D\qcd$
(i.e., to $\OX$-complexes whose homology sheaves are coherent).
Let  $\kappa=\kappa_Z$ be the canonical ringed-space map from the formal
completion~$X_{\</Z}$ to $X$, so that $\kappa_*$  and $\kappa^*$ are exact
functors \cite{GrD, p.\,422, (10.8.9)}. For $\F \in\nmb \A\qcd\>$,
following \cite{GrD, p.\,418, (10.8.2)} we denote by $\F_{\</Z}$ the
restriction of ~$\Lambda_Z\F$ to~$Z$.
{}From the map $\kappa^*\F\to \F_{\</Z}$ which is adjoint to
the natural map $\F\to \Lambda_Z\F=\kappa_*\F_{\</Z}\>$
we get a functorial map
$\kappa_*\kappa^*\F\to \nmb\kappa_*\F_{\</Z}\nmb=\nmb\Lambda_Z\F\>$;
and since $\kappa_*\kappa^*$ is exact, there results a functorial map
$$
\lambda_*^*\:\kappa_*\kappa^*\F\to \L\Lambda_Z\F\qquad (\F\in\D\qcd).
$$

\proclaim{Proposition (0.4.1)} The map\/ $\lambda_*^*$ is an isomorphism for
all\/~$\F\in\nmb\D_{\text{c}}$.
\endproclaim

\proof
The question being local, we may assume $X$ affine.
As indicated at the end of \S2, the functor ~$\L\Lambda_Z$ is bounded above
(i.e., ``way-out left'') and also bounded below
(i.e., ``way-out right''); and the same is clearly true
of~$\kappa_*\kappa^*\<$.
So by \cite{H, p.\,68, Prop.\,7.1} (dualized) we reduce to
where $\F$ is a single
finitely-generated free \hbox{$\OX$-module}, in~which case the assertion is
obvious since by \S1, $\L\Lambda_Z\P=\Lambda_Z\P$ for any quasi-coherent
flat complex~$\P$.
\endproof

Via the natural isomorphism
$
\kappa_*\RsH{X_{\</Z}}(\kappa^*\E,\>\kappa^*\F\>)\<\iso\<
\RsH{X}\<(\E,\>\kappa_*\kappa^*\F\>)
$
\cite{Sp, p.\,147, Prop.\,6.7\kern.6pt}, {\it the isomorphism in ~$(0.3)$ now
becomes, for\/} $\E\in\D,\ \F\in\Dc\!:$
$$
\qquad\RsH{X}\<(\R\vgp Z\E,\jog\F\>)\iso
 \kappa_*\RsH{X_{\<\</Z}}\<(\kappa^*\E,\>\kappa^* \F\>),
\tag "$(0.3)_{\text{\rm c}}$"
$$
{\it or\/}---by Lemma (0.4.2) below, and since as before
$\R\vgp Z\F\cong \R\vg Z\F\jog$:
$$
\RsH{X}\<(\R\vgp Z\E,\>\R \vgp Z\F\>)\iso
 \kappa_*\RsH{X_{\<\</Z}}\<(\kappa^*\E,\>\kappa^* \F\>).
\tag "$(0.3)_{\text{\rm c}}'$"
$$
\vskip -1\jot
\goodbreak
\noindent Explicitly, all these isomorphisms fit into a natural commutative
diagram:

\newcount\SaveCatCode

\SaveCatCode=\the\catcode`@
\catcode`@=11
\expandafter\ifx\csname graph\endcsname\relax
\alloc@4\box\chardef\insc@unt\graph\fi
\catcode`@=\number\SaveCatCode
\setbox\graph=\vtop{%
  \baselineskip=0pt \lineskip=0pt \lineskiplimit=0pt
  \vbox to0pt{\hbox{%
    \special{pn 8}%
    \special{pn 2}%
    \special{pa 0 249}%
    \special{pa 379 0}%
    \special{fp}%
    \special{pa 376 35}%
    \special{pa 379 0}%
    \special{fp}%
    \special{pa 346 -10}%
    \special{pa 379 0}%
    \special{fp}%
    \kern  0.380in
  }\vss}%
  \kern  0.250in
}




\SaveCatCode=\the\catcode`@
\catcode`@=11
\expandafter\ifx\csname graphb\endcsname\relax
\alloc@4\box\chardef\insc@unt\graphb\fi
\catcode`@=\number\SaveCatCode
\setbox\graphb=\vtop{%
  \baselineskip=0pt \lineskip=0pt \lineskiplimit=0pt
  \vbox to0pt{\hbox{%
    \special{pn 8}%
    \special{pn 2}%
    \special{pa 0 0}%
    \special{pa 379 249}%
    \special{fp}%
    \special{pa 346 260}%
    \special{pa 379 249}%
    \special{fp}%
    \special{pa 376 214}%
    \special{pa 379 249}%
    \special{fp}%
    \kern  0.380in
  }\vss}%
  \kern  0.261in
}
$$
\minCDarrowwidth{.3in}
\def\1{\RsH{X}\<(\R\vgp Z\E,\>\R \vgp Z\F\>)}
\def\2{\RsH{X}\<(\R\vgp Z\E,\>\F\>)}
\def\3{\RsH{X}\<(\E,\>\F\>)}
\def\4{\RsH{X}\<(\E,\>\L\Lambda_Z\F\>)}
\def\5{\RsH{X}\<(\E,\>\kappa_*\kappa^*\F\>)}
\def\6{\kappa_*\RsH{X_{\</Z}}(\kappa^*\E,\>\kappa^*\F\>)}
\CD
@.\2 @<\Iso<(0.4.2)<
    \underset{\UnderElement{\simeq}{\downarrow}{7.7ex}
    {\<\<(0.3)_{\text{c}}'}} \to\1 \\  \vspace{-4pt}
@. @V\kern-25mm{\box\graph}V\mkern-27.5mu\simeq\mkern14mu
 (0.3) V \\
\3 @>\hskip 1.8em>> \4 \\ \vspace{1.6pt}
@. @A\kern-25mm{\box\graphb}A\kern-4.5mm\simeq\kern2.5mm
 \text{via }\lambda_*^*A \\
\vspace{1pt}
@.\5 @<\Iso<\phantom{(0.4.2)}<\6
\endCD
$$
\proclaim{Lemma (0.4.2)} Let\/ $X$ be a scheme, let\/ $Z\subset X$ be a closed
subset, and let\/ $i\:(X\setminus Z)\hookrightarrow X$ be the inclusion.
Let\/ $\G\in\D(X)$ be exact off\/~$Z,$\ i.e., $i^*\G=\nmb0$.
Then for any\/ $\F\in\D(X)$ the natural map\/
$\RsH{}(\G,\>\R\vg Z\F\>)\to \RsH{}(\G,\>\F\>)$ is an isomorphism.
In particular, for any\/ $\E\in\D(X)$ there are natural isomorphisms
\vadjust{\kern-1\jot}
$$
\spreadlines{.5 \jot}
\gather
\RsH{}(\R\vg Z\E,\>\R\vg Z\F\>)\iso \RsH{}(\R\vg Z\E,\>\F\>),\\
\RsH{}(\R\vgp Z\E,\>\R\vg Z\F\>)\iso \RsH{}(\R\vgp Z\E,\>\F\>).
\endgather
$$
\endproclaim

\proof
If $\J$ is an injective K-injective resolution of~$\F$ \cite{Sp, p.\,138,
Thm.\,4.5} then $i^*\<\J$ is K-injective and the natural sequence $0\to\vg
Z\J\to\J\to i_*i^*\<\J\to 0$
is~exact; hence there is a natural triangle
$$
\R\vg Z\F\to\F\to\R\> i_*i^*\F@>\upl>>\R\vg Z\F\>[1]. \tag 0.4.2.1
$$
Apply the functor $\RsH{}(\G,\>-)$ to this triangle, and conclude
via the isomorphism
$\RsH{}(\G,\>\R\> i_*i^*\F\>)\cong \R\> i_*\RsH{}(i^*\G,\>i^*\F\>) =0$\ \;
\cite{Sp, p.\,147, Prop.\,6.7\kern .6pt}.
\endproof
\smallskip
{\bf(c)} (Local Duality). Let $A$ be a noetherian commutative ring
(so that any finite sequence in~$A$ is proregular), let $J$ be an
$A$-ideal, let $\hat A$ be the $J$-adic completion, and let
$\Gamma_{\!\!J}^{\phantom{.}}$
be the functor of $A$-modules described by
$$
\Gamma_{\!\!J}^{\phantom{.}}(M)\set
 \{\,x\in M\mid J^nx=0 \text{\quad\! for some }n>0\,\}.
$$
The derived $A$-module category $\D(\<A)$ has the full subcategory~$\Dc(\<A)$
consisting of those complexes whose homology modules are finitely generated.
 Arguing as in Remark~(b), one deduces from~$(0.3)_{\text{aff}}\jog$ the
 {\it duality isomorphism\/}
$$
\gathered
\R\Hom_{\mkern-1.5mu A}^\bullet(\R\Gamma_{\!\!J}^{\phantom{.}} E,
 \>\R\Gamma_{\!\!J}^{\phantom{.}}F\>)\iso
  \R\Hom_{\mkern-1.5mu A}^\bullet(E, \>F\otimes_A \<\hat A)\\
\hskip63pt\bigl(E\in\D(\<A),\; F\in\Dc(\<A)\bigr).\hskip-150pt
\endgathered\tag "$(0.3)_{\text{aff,\,c}}'$"
$$
(This is of course closely related to $(0.3)_{\text{c}}'\>$,
see Remark (d). For example,
when $J$~is a maximal ideal and $Z\set\{J\}\in X\set\text{Spec}(\<A)$,
just check out the germ of~$(0.3)_{\text{c}}'\>$~at the closed point
$J\in X$.)

\goodbreak
Now suppose that $E$ and $F$ are both in $\Dc(\<A)$, and one of the
following holds:
\smallskip
(1) $E\in\Dc^{\!\!\!\umi}(\<A)$ and $F\in\Dc^{\!\!\!\upl}(\<A)$;\kern 1pt
\footnote
{For any derived category $\D_{\<*}\>$,  $\D_{\<*}^{\<\upll}$
(resp.~$\D_{\<*}{}^{\mkern-10mu\umi}$) is the full subcategory whose objects
are the complexes~$C\in\D_{\<*}$
having bounded-below (resp.~bounded-above) homology\-, i.e., $H^n(C)=0$ for~
$n\ll0$ (resp.~$n\gg0$). $\D_{\<*}^{\<\upll}$
(resp.~$\D_{\<*}{}^{\mkern-10mu\umi}$) is isomorphic to
the derived category of the homotopy category of such ~$C$. This notation
differs from that in \cite{H}, where $C$ itself is assumed bounded.%
} or

(2) $F$ has finite injective dimension (i.e., $F$ is
$\D$-isomorphic to a bounded injective complex); or

(3) $E$ has finite projective dimension.
\smallskip
Then the natural map
$$
\nopagebreak
\R\Hom_{\mkern-1.5mu A}^\bullet(E,\>F\>)\otimes_{\!A}\<\hat
A\to\R\Hom_{\mkern-1.5mu A}^\bullet(E,\>F\otimes_A\<\hat A)
$$
is an isomorphism. To see this, reduce via ``way-out" reasoning
\cite{H, p.\,68}
to where $E$ is a bounded-above complex of finitely generated
projectives and $F$ is a single
finitely generated $A$-module. Similarly, $\Ext_{\<\<A}^n(E, \>F)\set
H^n\bigl(\R\Hom_{\<\<A}^\bullet(E,\>F\>)\bigr)$ is
 finitely generated.
Hence~$(0.3)_{\text{aff,\,c}}'$ yields homology isomorphisms
$$
\qquad\Ext_{\<\<A}^n(\R\Gamma_{\!\!J}^{\phantom{.}} E,
 \>\R\Gamma_{\!\!J}^{\phantom{.}}F\>)\iso
\Ext_{\<\<A}^n(E, \>F)\,\hat{}\qquad(n\in\Bbb Z).
$$

In particular, if $m$ is a maximal ideal and $D\in\Dc(\<A)$ is a
{\it dualizing complex\/} (which has, by definition finite injective
dimension), normalized so that $\R\Gamma_{\!\!\jog m}^{\phantom{.}}D$
is an injective hull~$I_m$
of the $A$-module $A/m$ \cite{H, p.\,284, Prop.\,7.3}, then there are
hyperhomology duality isomorphisms, generalizing \cite{H, p.\,280, Cor.\,6.5}:
$$
\Hom_{\<\<\hat A}(\Bbb H_{\>m}^{-n}E,\>I_m\>)\iso
 \text{Ext}_{\<\<A}^n(E, \>D)\,\hat{}
   \qquad\bigl(n\in\Bbb Z,\ E\in\Dc(\<A)\bigr).
$$
And since $\text{Ext}_{\<\<A}^n(E, D)\,\hat{}\ $ is a noetherian $\hat
A$-module
therefore $\Bbb H_{\>m}^{-n}E$ is artinian, and
Matlis dualization yields the Local Duality theorem of \cite{H, p.\,278}.
(One checks that the isomorphisms derived here
agree with those in~\cite{H}.)\looseness=-1

{\it More generally,} if $J$ is any $A$-ideal and
\lower.4ex\hbox{$\widehat{\phantom I}\>$} denotes $J$-adic completion
then with
$\kappa\:\text{Spf}\>(\hat A)=\nmb\widehat X\<\to@!@! X\set\text{Spec}(\<A)$
the canonical map, $\>U\set\nmb X\setminus\{m\}$,  and
$\>\E\set\nmb\w E$,
$\Cal D\set\nmb\w D$ the quasi-coherent $\OX$-complexes
generated by $E\/$ and~$D$,  there is a~triangle\looseness=-1
$$
\Hom_{\<\<A}^\bullet(\R\Gamma_{\!\!J}^{\phantom{.}}E,\> I_m\>)\to
 \R\Hom_{\<\<A}^\bullet(E,D)\otimes_{\<\<A}\<\<\hat A\to
 \R\Hom_{\widehat U}^\bullet(\kappa^*\E,\> \kappa^*\Cal D)@>\upl>>
$$
whose exact homology sequence looks like
$$
\dots\to\Hom_{\<\<\hat A}(\Bbb H_{\>J}^{-n}E,\>I_m\>)\to
 \text{Ext}_{\<\<A}^n(E,D)\,\hat{}\to
  \text{Ext}_{\widehat U}^n(\kappa^*\E,\> \kappa^*\Cal D)\to\dots
\tag 0.4.3
$$
The particular case when $A$ is Gorenstein of dimension~$d$---so that
$D\cong A[d\>\>]$---and $E=A$, is \cite{PS, p.\,107, Prop.\,(2.2)}.
See \S5.4 for details.

Incidentally, we have here a characterization
of $D\otimes_{\<\<A}\<\<\hat A$ as
$$
D\otimes_{\<\<A}\<\<\hat A\underset{(0.3)_{\text{aff,\,c}}'}\to\cong
\R\Hom_{\mkern-1.5mu A}^\bullet(\R\Gamma_{\!m}A,\mkern1mu \R\Gamma_{\!m}D)=
\R\Hom_{\mkern-1.5mu A}^\bullet(\R\Gamma_{\!m}A,\mkern1mu I_m)
\underset{(0.3)_{\text{aff}}}\to\cong
 \LL mI_m\>.
$$
Thus if $E^\bullet$ is an injective resolution of~$A$, so that
$\Hom_A(E^\bullet,\>I_m)$
is a flat resolution of~$I_m$ \cite{M, p.\,95, Lemma 1.4}, then
$D\otimes_{\<\<A}\<\<\hat A\cong \Hom_A(E^\bullet\<,\>I_m)\,\hat{}\,$.
\medskip
\eightpoint

{\bf(d)} Not surprisingly, but also not trivially,
$(0.3)_{\text{aff}}$ can be derived from~(0.3)---and vice-versa when
$\E\i\D\qcd(X)$---in brief as follows.

The functor $\GG X\set\Gamma(X,-)$ ($X\set\text{Spec}(\<A)$)
has an exact left adjoint, taking an
$A$-module~$M$ to its associated quasi-coherent $\OX$-module~$\w M$. Hence
$\GG X$ preserves K-injectivity, and there is a functorial isomorphism
$$
\R\Hom_{\mkern-1.5mu A}^\bullet(E,\>\R\GG X\G)\iso \R\Hom_{\<X}^\bullet(\w
E,\>\G)
\qquad\bigl(E\in\D(\<A),\; \G\in\D(X)\bigr).
$$

Next, if $\G$ is any $\OX$-complex of $\GG X$-acyclics
(i.e., the natural map $\GG X\G^n\to\R\GG X\G^n$ is an isomorphism
for all~$n$), then  $\GG X\G\to\R\GG X\G$ is an isomorphism.
(This is well-known\vadjust{\nopagebreak} if $\G$ is bounded below; and in the
general case can be deduced from \cite{BN, \S5} or found explicitly in
\cite{L, (3.9.3.5)}.) So for any $A$-complex~$F$ there are natural isomorphisms
$F\iso\GG X\w F\iso\R\GG X\w F$, and hence
$$
\R\Hom_{\mkern-1.5mu A}^\bullet(E,\>F)\iso \R\Hom_{\<X}^\bullet(\w E,\>\w F)
\qquad\bigl(E,F\in\D(\<A)\bigr). \tag 0.4.4
$$
There are also natural isomorphisms
$$
\R\vg Z\w E\iso \w{\R\GG{\>\>\bt}\<\<E}, \qquad\
\LL\bt F\iso \R\GG X\LL Z\w F.\tag 0.4.5
$$
The first obtains via Koszul complexes, see (3.2.3). For the second, we
may assume ~$F$ flat and K-flat, in which case we are saying that
$\Lambda_\bt F=\GG X\Lambda_Z\w F\to \R\GG X\Lambda_Z\w F\>$ is an isomorphism,
which as above reduces to where $F$ is a single flat $A$-module,
and then follows from \cite{EGA, p.\,68, (13.3.1)}.

Thus there are natural isomorphisms
$$
\spreadlines{1\jot}
\gathered
\R\Hom_{\mkern-1.5mu A}^\bullet(E,\mkern1mu\L\Lambda_\bt F\>) \iso
 \R\Hom_{\mkern-1.5mu A}^\bullet(E,\mkern1mu\R\GG X\LL Z\w F\>)
\iso \R\Hom_X^\bullet(\w E,\LL Z\w F\>), \\
\R\Hom_{\mkern-1.5mu A}^\bullet\bigl(\R\Gamma_{\!\bt}\<\<E,\>F\>)\iso
\R\Hom_{\mkern-1.5mu A}^\bullet\bigl(\w{\R\GG{\>\>\bt}E},\>\w F\>)\iso
 \R\Hom_X^\bullet\bigl(\R\vg Z\w E,\>\w F\>).
\endgathered\tag "\tt (\#)"
$$
Hence (0.3) implies $(0.3)_{\text{aff}}$. Conversely, (0.3)(bis) (with
$\E\in\D\qcd(X))$ follows from $(0.3)_{\text{aff}}$. Indeed, it suffices to see
that the maps $\lambda'$ and $\gamma'$ are made into isomorphisms by
the functor~$\R\GG {\>U}$ for any affine open $U\subset X$.
Moreover, we may assume that the
complexes $\E$ and~$\F$ are quasi-coherent (see \S1). Then
{\tt (\#)} provides what we need.
\par\tenpoint
\medskip
\goodbreak

\subheading{1. Left-derivability of the completion functor}
Let $X$ be a quasi-compact separated scheme
and let $Z\subset X$ be a closed subscheme.
We show in this section that
{\it the completion functor\/ $\Lambda_Z\:\A\qcd(X)\to\A(X)$
of\/~$(0.2)$ has a left-derived functor\/
$\L\Lambda_Z\:\D\qcd(X)\to\D(X)$.}

\proclaim{Proposition (1.1)} On a quasi-compact separated scheme\/~$X,$\
every $\E\in\nmb\D\qcd(X)$ is isomorphic  to a quasi-coherent
K-flat complex~$\P_{\!\E}$.
\endproclaim

The {\it proof\/} will be given below, in~(1.2).
\smallskip
If $\P\in\D(X)$ is a K-flat
exact {\it  quasi-coherent\/} complex, then\/ $\Lambda_Z(\P)$ is exact.
\hbox{Indeed,} all the complexes
$\P_n\set (\O/\I^{\>n})\jog\otimes \jog \P\ \;(n>0)$
are exact \cite{Sp, p.\,140, Prop.\,5.7},
and hence the same is true after taking global sections
over any affine open subset~$U$ of~$X\<$.
Also, the natural map of complexes
$\Gamma(U, \P_{n+1})\to\nmb\Gamma(U, \P_n)$ is surjective for every~$n$.
So by \cite{EGA, p.\,66, (13.2.3)}, the complex
$$
\Gamma(U, \Lambda_Z(\P))=\inlm\,\Gamma(U, \P_n)
$$
is exact, whence the assertion.

Consequently (see \cite{H, p.\,53}, where condition~1
for the triangulated subcategory~$L$ whose objects are
all the quasi-coherent K-flat complexes can be replaced
by the weaker condition in our Proposition~(1.1)), after choosing one
$\P_{\!\E}$ for each $\E$  we have a left-derived functor $\L\Lambda_Z$ with
$\L\Lambda_Z(\E)\set\Lambda_Z(\P_{\!\E})$.
For simplicity we take $\P_{\!\E}=\E$ whenever $\E$~itself is quasi-coherent
and
K-flat, so then $\L\Lambda_Z(\E)=\Lambda_Z(\E)$.

\eightpoint
\nextpart{1.2} Here is the {\it proof of Proposition\/~$(1.1)$.}
It uses a simple-minded version of some simplicial techniques
found e.g., in \cite{Ki, \S2}. We will recall as much as is needed.

Let $\Cal U=(U_\alpha)_{1\le \alpha\le n}$\vadjust{\kern.5pt}
be an affine open cover
of the quasi-compact separated scheme~$(X,\OX\<)$. Denote the set of
subsets of $\{1,2,\dots,n\}$ by~$\fP_n$. For $i\in \fP_n$, set
$$
U_i\set \bigcap_{\alpha\in i}U_\alpha\>,\qquad \O_i\set\O_{U_i}=\OX|_{U_i}\>.
$$
(In particular, $U_\phi=X$.) For $i\supset j$ in~$\fP_n$, let
$
\lambda_{ij}\:U_i\hookrightarrow U_j
$
be the inclusion map.
A $\Cal U$-module is, by definition, a family
$\F=(\F_i)_{i\in{\ssize\frak P_n}}$
where $\F_i$ is an $\O_i$-module, together with
a family of $\O_j$-homomorphisms

$$\varphi_{jk}\:\lambda_{jk}^*\F_k\to\F_j\qquad(j\supset k)
$$
such that $\varphi_{jj}$ is the identity map of~$\F_j$, and whenever $i\supset
j\supset k$ we have
$\varphi_{ik}=\nmb\varphi_{ij}\smcirc\nmb\bigl(\varphi_{jk}|_{U_i}\bigr)$,
i.e.,
$\varphi_{ik}$ factors as
$$
\lambda_{ik}^*\F_k=\lambda_{ij}^*\lambda_{jk}^*\F_k
@>\lambda_{ij}^*(\varphi_{jk})>>\lambda_{ij}^*\F_j
@>\varphi_{ij}>>\F_i\jog.
$$
We say the $\Cal U$-module $\F$ is {\it quasi-coherent\/}
(resp.~{\it flat\/}, resp.~\dots) if each one of the $\O_i$-modules~$\F_i$
is such.

The $\Cal U$-modules together with their morphisms (defined in the obvious
manner) form an abelian category with $\drlm$ and $\inlm$. For example, given
a direct system $(\F^\rho)_{\rho\in R}$ of $\Cal U$-modules,
set $\F_i\set\drlm\F_i^\rho\ (i\in \fP_n)$,
define $\varphi_{ij}$ ($i\supset\nmb j\>$) to be the adjoint of
the natural composed map
$$
\F_j=\drlm \F_j^\rho @>\text{via }\psi_{ij}^\rho>>
\drlm\lambda_{ij*}\F_i^\rho\lra \lambda_{ij*}\F_i
$$
where $\psi_{ij}^\rho\:\F_j^\rho\to \lambda_{ij*}\F_i^\rho$ is adjoint to
$\varphi_{ij}^\rho\:\lambda_{ij}^*\F_j^\rho\to \F_i^\rho$; and check
that \smash{$\F\set(\F_i@,@,,\jog\varphi_{ij})=\nmb\drlm\<\F^\rho$} in the
category of
$\Cal U$-modules.

\proclaim{Lemma(1.2.1)}
Any quasi-coherent\/ $\Cal U$-module\/ $\F$ is a homomorphic image of a\/
flat quasi-coherent $\Cal U$-module.
\endproclaim

\demo{Proof}
For each $i$ we can find an epimorphism of quasi-coherent $\O_i$-modules
$\Q_i\twoheadrightarrow \F_i$ with $\Q_i$ flat. Set
$\P_i\set\oplus_{i\supset j}\jog\lambda_{ij}^*\Q_j$. Map $\P_i$ surjectively
to~
$\F_i$ via the family of composed maps
$$
\lambda_{ij}^*\Q_j\lra \lambda_{ij}^*\F_j@>\varphi_{ij}>>\F_i.
$$
Let
$$
\varphi_{ki}'\:\lambda_{ki}^*\P_i=\oplus_{i\supset j}\jog\lambda_{kj}^*\Q_j
\lra \oplus_{k\supset j}\jog\lambda_{kj}^*\Q_j =\P_k
$$
be the natural map. Then $\P\set(\P_i@,@,,\jog\varphi_{ij}')$ is a flat
$\Cal U$-module, and the maps $\P_i\to \F_i$ constitute an epimorphism of
$\Cal U$-modules.
\endproof

The tensor product of two $\Cal U$-modules is defined in the obvious way. A
complex of $\Cal U$-modules is {\it K-flat} if its tensor product with any
exact
complex is again exact.

\proclaim{Corollary(1.2.2)} {\rm(Cf.~\cite{Ki, p.\,303, Satz 2.2.})}
Any complex of quasi-coherent\/ $\Cal U$-modules is the target of a
quasi-isomorphism from a K-flat complex of quasi-coherent\/ $\Cal U$-modules.
\endproclaim

\demo{Proof} (Sketch.) Any bounded-above complex of flat $\Cal U$-modules is
K-flat, so the assertion for bounded-above complexes
follows from Lemma~(1.2.1)
(see \cite{H, p.\,42, 4.6, 1) (dualized)}). In the general case,
express an arbitrary complex as the \smash{$\drlm$}\vadjust{\kern.7pt}
of its truncations,
and then use the \smash{$\drlm$}\vadjust{\kern.7pt}
of a suitable direct system of K-flat resolutions
of these truncations. (Clearly, \smash{$\drlm$}\vadjust{\kern.5pt} preserves
K-flatness.
For more details, see \cite{Sp, p.\,132, Lemma\,3.3}
or \cite{L, (2.5.5)}.)
\endproof
The {\it \v Cech functor} $\check C^\bullet$ from $\Cal U$-complexes
(i.e., complexes of $\Cal U$-modules) to $\OX$-complexes is defined
as follows:

Let $|i|$ be the cardinality of $i\in \fP_n$, and let
$\lambda_i\set\lambda_{i\phi}$ be the inclusion map $U_i\hookrightarrow X\<$.
For any $\Cal U$-module ~$\F\<$, set
$$
\spreadlines{1.2\jot}
\alignat2
\C^m(\F\>)&\set \bigoplus_{|i|=m+1} \lambda_{i*}\F_i\qquad&&0\le m<n\\
&\set 0 &&\text{otherwise.}
\endalignat
$$
Whenever $j$ is obtained from $k=\{k_0<k_1<\dots<k_m\}\in \fP_n$ by removing a
single element, say ~$k_a$, we set $\epsilon_{kj}\set(-1)^a$. The boundary map
$\delta^m\:\C^m(\F\>)\to\nmb\C^{m+1}(\F\>)$ is specified by the family of maps
$$
\delta_{kj}^m\: \lambda_{j*}\F_j\to \lambda_{k*}\F_k
$$
with $\delta_{kj}^m$ the natural composition
$$
\nopagebreak
\lambda_{j*}\F_j
 \lra\lambda_{j*}\lambda_{kj*}\lambda_{kj}^*\F_j
  =\lambda_{k*}\lambda_{kj}^*\F_j
   @>\lambda_{k*}(\epsilon_{kj}\varphi_{kj})>>
    \lambda_{k*}\F_k
$$
if $j\subset k$, and $\delta_{kj}^m=0$ otherwise. Then
$\delta^{m+1}\smcirc@,@,\delta^m=0$ for all~$m$, and so we have
a functor $\C^\bullet$ from $\Cal U$-modules to $\OX$-complexes.
For any \hbox{$\Cal U$-complex~$\F^\bullet\<$}, $\C^\bullet(\F^\bullet)$ is
defined to
be the total complex\vadjust{\kern.7pt} associated to the double
complex~$\C^p(\F^q)$.

\remark{Remarks} (a) If $\G$ is an $\OX$-module and $\G'$ is the
$\Cal U$-module such that $\G_i'\set\nmb\lambda_i^*\G$ and $\varphi_{ij}$ is
the identity map of $\G'_i=\lambda_{ij}^*\G_j$ for all $i\supset j$, then
$\C^\bullet(\G')$ is the usual \v Cech resolution of ~$\G$
\cite{Go, p.\,206, Thm.\,5.2.1}.

(b) Since all the maps $\lambda_i$ are {\it affine\/}
($X$ being separated) and {\it flat},
therefore $\C^\bullet$ takes flat quasi-coherent $\Cal U$-complexes
to flat quasi-coherent $\OX$-complexes. Moreover, $\C^\bullet$
commutes with \smash{$\drlm$}\vadjust{\kern.7pt}.
(We need this only for quasi-coherent complexes,
for which the proof is straightforward;
but it also holds for arbitrary complexes, \cite{Ke, \S2}.)
\endremark

\proclaim{Lemma (1.2.3)}
The functor $\C^\bullet$ takes quasi-isomorphisms between\/ {\rm
quasi-coherent} complexes to quasi-isomorphisms.
\endproclaim

\demo{Proof} \kern-.7pt One checks that $\C^\bullet$ commutes with
degree-shifting:
$\C^\bullet(\F^\bullet[1])=\C^\bullet(\F^\bullet)[1]$; and that
$\C^\bullet$ preserves mapping cones. Since quasi-isomorphisms are just those
maps whose cones are exact, it suffices to show that $\C^\bullet$ takes
exact quasi-coherent $\Cal U$-complexes~$\F^\bullet$ to exact $\OX$-complexes.
But since the  maps $\lambda_i$ are affine, each row $\C^p(\F^\bullet)$ of
the double complex~$\C^p(\F^q)$ is exact, and all but finitely many rows
vanish, whence the conclusion. \endproof
\smallskip

Now by \cite{BN, p.\, 230, Corollary 5.5}, any $\E\in\D\qcd(X)$ is isomorphic
to a quasi-coherent complex; so to prove (1.1) we may as well assume
that $\E$ itself is quasi-coherent.
Define the $\Cal U$-complex $\E'$ as in
remark~(a) and let $\P\to \E'\>$ be a quasi-isomorphism of
quasi-coherent $\Cal U$-complexes
with $\P$ a \smash{$\drlm$}\vadjust{\kern.7pt} of bounded-above flat complexes,
see proof of Corollary~(1.2.2). Lemma~(1.2.3) provides
a quasi-isomorphism $\P_\E\set\C^\bullet(\P)\to\C^\bullet(\E')$; and there is a
natural  quasi-isomorphism $\E\to\C^\bullet(\E')$ (remark~(a)), so that
$\E$ is isomorphic in~$\D(X)$ to $\P_\E$.
Moreover, $\P_\E$ is a \smash{$\drlm$}
of bounded-above quasi-coherent flat $\OX$-complexes
(remark~(b)), and hence is quasi-coherent and K-flat. This proves
Proposition~(1.1). \qed
\medskip

For completeness, and for later use, we present a slightly more
elaborate version of the just-quoted Corollary~5.5 in \cite{BN, p.\,230}.
Recall from Remark (0.4)(a) the definition of {\it quasi-coherator.}

\proclaim{Proposition (1.3)}
Let $X$ be a quasi-compact separated scheme. Then the natural functor
$$
{\boldkey j}_X\:\D(\A\qcd(X)) \to \D\qcd(X)
$$
is an equivalence of categories, having as quasi-inverse the derived
quasi-coherator\/~$\Rq X$.
\endproclaim
\proclaim{Corollary (1.3.1)}
In the category\/~
$\text{\eighteurm C}\qcd(X)$ of quasi-coherent $\O_X$-complexes,
every object has a K-injective resolution.
\endproclaim
\proof
The Proposition asserts that the natural maps $\E\to\Rq X\boldkey j_X\E$
$\bigl(\E\in\D(\A\qcd(X))\bigr)$ and $\boldkey j_X\Rq X\F\to\nmb\F$
$\bigl(\F\in\nmb\D\qcd(X)\bigr)$ are isomorphisms. The Corollary results:
since $Q_X$ has an  exact left adjoint therefore $Q_X$ takes
K-injective $\O_X$-complexes to complexes which are K-injective in~
$\text{\eighteurm C}\qcd(X)$, so if $\E\iso\Rq X\boldkey j_X\E\>$ and
if $\E\to I_\E$ is a quasi-isomorphism with  $I_\E$
a \hbox{K-injective} $\O_X$-complex \cite{Sp, p.\,134, 3.13},
then the resulting map
$\E\to Q_XI_\E$ is still a quasi-isomorphism, and thus $\E$ has a K-injective
resolution in~$\text{\eighteurm C}\qcd(X)$.

We will show that
the functor $\Rq X|_{\D\qcd(X)}$ is {\it bounded-above,} i.e., there
is an integer $d\ge0$ such
that for any $\F\in\D\qcd(X)$ and $q\in\text{\eightBbb Z}$, if
$H^p(\F)=0$ for all $p\ge q$ then $H^p(\Rq X\F)=0$ for all $p\ge q+d$.
Then by the way-out Lemma \cite{H, p.\,68} it suffices to prove the above
isomorphism assertions when $\E$ and $\F$ are single quasi-coherent sheaves,
and this case is dealt with in \cite{I, p.\,189, Prop.\,3.5}. (It follows
then from $\boldkey j_X\Rq X\F\iso\F$ that we can take $d=0$.)

We proceed by induction on $n(X)$, the least among all
natural numbers~$n$ such that $X$ can be covered by $n$ affine open subschemes.
If $n(X)=1$, i.e., $X$ is affine,
then for any $\F\in\D\qcd(X)$, $\Rq X(\F\>)$ is the sheafification
of the complex $\R\GG X(\F\>)\set\R\Gamma(X,\F\>)$;
so to show boundedness we can replace
$\Rq X$ by $\R\GG X$. For a K-injective resolution~$I$ of~$\F\in\D\qcd(X)$,
use a ``special" inverse limit of injective
resolutions~$I_q$ of the truncations $\tau\>_{{\sssize\ge} -q}(F\>)$, as in
\cite{Sp, p.\,134, 3.13}. If $C_q$ is the kernel of the split surjection
$I_q\to I_{q-1}$, then $C_q[-q]$ is an injective resolution of
the quasi-coherent $\OX$-module~$H^{-q}(\F)$, and hence
$H^p\GG X(C_q)=0$ for $p> -q$. Applying \cite{Sp, p.126, Lemma},
one finds then that for $p\ge -q$  the natural map
$H^p\GG X(I)\to H^p\GG X(I_q)$ is an isomorphism; and so if
$\tau\>_{{\sssize\ge} -q}(\F\>)=0$, then $H^p\GG X(I)=\nmb0$.
Thus $\R\GG X|_{\D\qcd(X)}$ is indeed bounded above (with $d=0$).

Now suppose that $n\set n(X)>1$, and let $X=X_1\cup\dots\cup X_n$
be an affine open cover. Set $U\set X_1$, $V\set X_2\cup\dots\cup X_n\>$,
$W\set U\cap V$, and let $u\:U\hookrightarrow X$, $v\:V\hookrightarrow X$,
$w\:W\hookrightarrow X$ be the inclusions. Note that $n(U)=1$, $n(V)=n-1$,
and $n(W)\le n-1$ ($X$ separated $\Rightarrow X_1\cap X_i$ affine).

By the inductive hypothesis, $\E\iso\Rq V\boldkey j_V\E\>$ for any
$\E\in \text{\eighteurm C}\qcd(V)$. Hence, as above, $\E$ has a K-injective
resolution in~$\text{\eighteurm C}\qcd(V)$, the functor
$v_*\qcu\:\D(\A\qcd(V))\to\D(\A\qcd(X))$
has a right-derived functor~$\R v_*\qcu\<$, and there is a functorial
isomorphism $\R(v_*\qcu Q_V)\iso\nmb\R v_*\qcu@,\Rq V\>$.
Since the left adjoint $v^*$ of $v_*$ is exact,
therefore $v_*$ preserves
K-injectivity of complexes, and so there is a functorial isomorphism
$\R(Q_Xv_*)\iso\Rq X\R v_*$; and furthermore it is easily seen,
via adjointness of $v^*$ and $v_*$, that $Q_Xv_*=v_*\qcu Q_V$.
Thus we have a functorial isomorphism
$$
\nopagebreak
\Rq X\R v_*\iso \R(Q_Xv_*)=\R(v_*\qcu Q_V)\iso\R v_*\qcu@, \Rq V.
$$

Similar remarks apply to $u$ and $w$.

Now we can apply~$\Rq X$ to the
Mayer-Vietoris triangle
$$
\F\to\R u_*u^*\F\oplus \R \>v_*v^*\F\to \R\> w_*w^*\F \to \F[1]
$$
to get the triangle
$$
\Rq X\F\to \R u_*\qcu \Rq U u^*\F \oplus \R v_*\qcu \Rq V v^*\F \to
 \R w_*\qcu \Rq W w^*\F \to \Rq X\F[1].
$$
So it's enough to show: {\sl if\/ $V$ is any
quasi-compact open subset of\/~$X$ with\/ $n(V)<n(X),$\ and\/
$v\:V\hookrightarrow X$  is the inclusion,
then the functor\/ $\R v_*\qcu$ is bounded above.}
(This derived functor exists, as before, by the induction hypothesis.)

We induct on~$n(V)$, the case $n(V)=\nmb 1$ being trivial,
since then the map~$v$ is affine and the
functor $v_*\qcu\:\A\qcd(V)\to \A\qcd(X)$ is exact.
Suppose then that $n\set n(V)>1$.
$V$~has an open cover
$V=V_1\cup V_2$ with $n(V_1)=1$, $n(V_2)=n-1$, and $n(V_1\cap V_2)\le n-1$.
Let $i_1\:V_1\hookrightarrow V$, $i_2\:V_2\hookrightarrow V$, and
$i_{12}\:V_{12}= V_1\cap V_2\hookrightarrow V$ be the inclusions.
By \cite{L, (3.9.2)} (which uses techniques from \cite{Sp}
like those in the
above discussion of the case $n(X)=1$), $\R i_{1*}\D\qcd(V_1)\subset
\D\qcd(V)$,
and similarly for $i_2$ and~$i_{12}$.
Since $n(V\!_s)<n(X)\ (s=1,2, \text{ or }12)$,
we may assume that
${\boldkey j}_{V\!_s}\:\D(\A\qcd(V\!_s)) \to\nmb \D\qcd(V\!_s)$ is an
equivalence of categories with quasi-inverse~$\Rq {V\!_s}@!@!,@,$\ so that we
have isomorphisms
$$
\Rq V\R i_{s*}i_s^*\boldkey j_V\E \cong
 \R i_{s*}\qcu\Rq{V\!_s}{\boldkey j}_{V\!_s}i_s^*\E \cong
  \R i_{s*}\qcu i_s^*\E \qquad\bigl(\E\in\D(\A\qcd(V)\bigr).
$$
Similarly, $\Rq V\boldkey j_V\E\cong\E$. Hence application of~$\Rq V$
to the Mayer-Vietoris triangle
$$
\boldkey j_V\E\to
 \R i_{1*}i_1^*\boldkey j_V\E \oplus \R i_{2*}i_2^*\boldkey j_V\E \to
  \R i_{12*}i_{12}^*\boldkey j_V\E \to\boldkey j_V\E[1]
$$
gives rise to a triangle
$$
\E\to
 \R i@,@,_{1*}\qcu i_1^*\E \oplus \R i@,@,_{2*}\qcu i_2^*\E \to
  \R i@,@,_{12*}\qcu i_{12}^*\E \to\E[1].
$$
Since $i@,@,_{s*}\qcu$ has an exact left adjoint~$i_s^*$, therefore
$i@,@,_{s*}\qcu$ preserves K-injectivity, and consequently
$\R v_*\qcu\R i@,@,_{s*}\qcu=\nmb\R(vi_s)_*\qcu.$
So we can
apply $\R v_*\qcu$ to the preceding triangle and use the induction hypothesis
to see that $\R v_*\qcu\E$ is one vertex of a triangle whose other
two vertices are obtained by applying bounded-above functors to~$\E\<$,
whence the conclusion.
\endproof
\par\tenpoint

\subheading{2. Proof of Theorem~(0.3)---outline}
We first define bifunctorial maps
$$
\spreadlines{1\jot}
\aligned
\psi\:\E\Otimes\R\vg Z\F\to \R\vg Z(\E\Otimes\F\>) \\
\psi'\:\E\Otimes\R\vgp Z\F\to \R\vgp Z(\E\Otimes\F\>)
\endaligned
 \qquad\bigl(\E,\F \in\D(X)\bigr)\tag 2.1
$$
(where \smash{$\Otimes$} denotes derived tensor product.)
To do so, we may assume that $\E$ is K-flat and $\F$ is
K-injective, and choose a quasi-isomorphism $\E\otimes\F\to\J$ with
$\J$ K-injective. The obvious composed map of complexes
$
\E\otimes \vg Z\F\to \E\otimes\F\to\J
$
has image in $\vg Z \J$, and so we can define $\psi$ to be the
resulting composition in $\D(X)$:
$$
\botsmash{\E\Otimes\R\vg Z\F}\cong \E\otimes \vg Z\F\to \vg Z\J
 \cong\R\vg Z(\E\Otimes\F\>).
$$
The map $\psi'$ is defined similarly, {\it mutatis mutandis.}

Under the  hypotheses of~Theorem~(0.3), assertion~(i) in Cor.\,(3.2.5)
(resp.~(3.1.5)) gives that {\it $\psi$ is an isomorphism if\/
$\E$ and $\F$ are both in\/}~$\D\qcd(X)$ (resp.~ {\it $\psi'$ is an isomorphism
for all\/} $\E,\F\>)$.
\footnote
{The ring-theoretic avatar of this result
is closely related to results of
Matlis \cite{M, p.\,114, Thm.\,3.7}, \cite{M2, p.\,83, Thm.\,10},
and Strebel \cite{St, p.\,94, 5.8}.%
}

In view of the canonical isomorphism $\R\vgp Z\OX\iso\R\vg Z\OX$
(Cor.\,(3.2.4)) and of \cite{Sp, p.\,147, Prop.\,6.6}, we have then
natural isomorphisms
$$
\align
\RsH{}\bigl(\R\vgp Z\E,\> \F\>\bigr)&\iso
 \RsH{}\bigl(\E\Otimes \R\vg Z\OX,\> \F\>\bigr)\\
&\iso
  \RsH{}\bigl(\E,\mkern1mu \RsH{}(\R\vg Z\OX,\>\F\>)\bigr).
\endalign
$$

It remains to find a natural isomorphism
$$
\RsH{}(\R\vg Z\OX,\mkern1mu \F\>)\iso
 \L\Lambda_Z \F\qquad\bigl(\F\in\D\qcd(X)\bigr).
$$
To get this we define below a natural map
$\Phi\:\L\Lambda_Z \F\to\RsH{}(\R\vg Z\OX,\> \F\>)$,
and, after reducing to where $X$ is affine and
$\F$ is a single flat quasi-coherent
$\OX$-module,  prove in \S4 that $\Phi$ is an isomorphism
by constructing~$\Phi^{-1}$
via the representability of~$\R\vg Z\OX$
as a limit of Koszul complexes.
\footnote
{That proof turns out, at least for us, to be surprisingly difficult.
(A much shorter proof for the case $\F=\OX$ over smooth algebraic
${\text{\eightBbb C}}$-varieties is given by Mebkhout in \cite{Me, p.\,97}.
We could not follow all the details of his argument.)
The ring-theoretic avatar of the isomorphism~$\Phi$ underlies the
duality theorem of Strebel \cite{St, p.\,94, 5.9} and Matlis
\cite{M2, p.\,89, Thm.\,20},  and the more general results of
Greenlees and May \cite{GM, p.\,449, Prop.\,3.1 and p.\,447, Thm.\,2.5}.%
}
\medskip
Assuming  $X$ to be quasi-compact and separated,
so that $\L\Lambda_Z$ exists, let us then define $\Phi$. Let $\I$ be a
finite-type quasi-coherent $\OX$-ideal such that $Z=\text{Supp}(\OX/\I)$
(see Introduction).
For any $\OX$-complexes $\P$, $\Q$, $\Cal R$, the natural map
$$
\bigl(\P\otimes\Q\bigr)\otimes\bigl(\sHomb(\Q,\jog\Cal R)\bigr)\cong
\P\otimes\bigl(\Q\otimes\sHomb(\Q,\jog\Cal R)\bigr) \to \P\otimes\Cal R
$$
induces (via $\otimes$--$\sHom$ adjunction) a functorial map
$$
\P\otimes\Q\to\sHomb\bigl(\sHomb(\Q,\Cal R),\jog \P\otimes\Cal R\bigr)\jog.
$$
Letting $\Q$ run through the inverse system $\OX/\I^{\>n}\ \jog (n>0)$
one gets a natural map
$$
\align
\Lambda_Z(\P)=\inlm(\P\otimes\OX/\I^{\>n})
 &\to @,\inlm\sHomb\bigl(\sHomb(\OX/\I^{\>n},\Cal R),\jog \P\otimes\Cal
R\bigr)\\
&\cong\sHomb\bigl(\drlm\sHomb(\OX/\I^{\>n},\Cal R),\jog \P\otimes\Cal R\bigr)\\
&\cong\sHomb\bigl(\vgp Z\Cal R,\jog \P\otimes\Cal R\bigr)\jog.
\endalign
$$
For $\F\in\D\qcd(X)$, $\G\in\D(X)$,
taking $\P$ to be $\P_\F$
(Proposition (1.1)) and $\Cal R$ to be a K-injective
resolution of~$\G$ one gets a composed derived-category map
$$
\nopagebreak
\aligned
\Phi(\F\<,\G)\:
\L\Lambda_Z\F\cong \Lambda_Z\P
&\to \sHomb\bigl(\vgp Z\Cal R,\jog \P\otimes\Cal R\bigr)\\
&\to \RsH{}\bigl(\vgp Z\Cal R,\jog \P\otimes\Cal R\bigr)\\
&\cong \;\RsH{}\bigl(\R\vgp Z\G,\jog \F\Otimes\G\bigr)\jog,
\endaligned\tag 2.2
$$
which one checks to be independent of the choice of~$\P$ and~$\Cal R$.

As indicated above we want to show
that $\Phi(\F\<,\OX\<)$ is an isomorphism.
The question is readily seen to be local on~$X$,
\footnote
{Using the exact functor ``extension by zero," one shows that restriction
to an open~$U\subset\nmb X$ takes any K-injective (resp\. K-flat)
$\OX$-complex to a K-injective (resp\. K-flat) $\O_U$-complex.\looseness=-1}
so we may assume $X$ to be affine. The idea is then to apply
way-out reasoning \cite{H, p.\,69, (iii)} to reduce to where
$\F$ is a single {\it flat\/} quasi-coherent $\OX$-module, which case is
disposed of in~Prop.\,(4.2).

But to use {\it loc.\,cit.,} we need the functors
$\H_Z\set\RsH{}(\R\vg Z\OX,\jog-)$ and $\L\Lambda_Z$
from\/ $\D\qcd(X)$ to\/ $\D(X)$ to be bounded
above (=\medspace``way-out left'') and also bounded below
(=\medspace``way-out right''). Boundedness of~$\H_Z$ is shown in
Lemma~(4.3). That $\L\Lambda_Z(-)$ is bounded above is clear, since
$X$ is now affine and so  if $\E\in\D\qcd(X)$ is such that
$H^i(\E)=0$ for all $i>i_0$ then there is a flat $P_{\E}$ as in~(1.1)
vanishing in all degrees $>i_0\>$.
Now by \cite{H, p.\,69, (ii), (iv)} (dualized),
the case where $\F$ is a  flat quasi-coherent $\OX$-module (Prop.\,(4.2))
implies that $\Phi(\F,\OX\<)\:\H_Z\F\to\L\Lambda_Z\F$
is an isomorphism for all $\F\in\D\qcd^-(X)$.
Knowing that, and the fact that $\H_Z$ is bounded below, we can
conclude, by \cite{H, p.\, 68, Example~1} (dualized, with $P$ the class
of quasi-coherent flat $\OX$-modules), that $\L\Lambda_Z$ is bounded
below.  (See also \cite{GM, p.\,445, Thm.\,1.9,~(iv)}.)
\smallskip
For the last assertion of Theorem~(0.3), it suffices to verify the
commutativity
of the following diagram, where $\E$ may be taken to be K-flat, and
as above,  $\P=\P_\F\>$. This verification is straightforward (though
not entirely effortless) and so will be left to the reader.
$$
\nopagebreak
\def\1{\RsH{}\bigl(\E,\jog \P\>\bigr)}
\def\2{\RsH{}\bigl(\E,\jog \RsH{}\>(\OX,\>\P\>)\bigr)}
\def\3{\RsH{}\bigl(\E\otimes\OX,\jog \P\>\bigr)}
\def\4{\RsH{}\bigl(\E,\jog \F\>\bigr)}
\def\5{\RsH{}\bigl(\E,\jog \Lambda_Z\P\bigr)}
\def\6{\RsH{}\bigl(\E,\jog \RsH{}\>(\R\vgp Z\OX,\>\P\>)\bigr)}
\def\7{\RsH{}\bigl(\E\otimes\R\vgp Z\OX,\jog \P\>\bigr)}
\def\8{\RsH{}\bigl(\R\vgp Z\E,\jog \F\>\bigr)}
\CD
\1 @>\text{via }\lambda>>\5 \\
@| @V\text{via} V \Phi(\P\!,\>\OX) V  \\
\2 @>>> \6  \\
@| @VV\simeq V  \\
\3 @>>> \7 \\
@V\simeq VV @V(2.1)V\simeq V  \\
\4 @>>\text{via }\gamma> \8
\endCD
$$

This completes the outline of the proof of Theorem~(0.3).

\medskip

\subheading{3. Proregular embeddings}
In this section we explore the basic condition of {\it proregularity,}
as defined in (3.0.1). This definition, taken from \cite{GM, p.\,445}, seems
unmotivated at first sight; but as mentioned in the Introduction, it is
precisely what is needed to make local cohomology on quite general schemes
behave as it does on noetherian schemes (where {\it every\/} closed subscheme
is proregularly embedded), for example with respect to Koszul complexes.
What this amounts to basically is an elaboration of \cite{Gr,
Expos\'e~II\kern.6pt} in the language of derived categories of sheaves.
\footnote
{More generally,
to do the same for [{\it ibid.,} Expos\'e~VI\kern.6pt], replace $\O$ in what
follows by an
\hbox{$\O$-module} $\Cal M$, $\P$ by $\Cal M\otimes\P$ ($\P$ flat), $\J$ by
$\H\text{\eighteurm om}\>\>(\Cal M, \J)$ ($\J$ injective), and the functor
$\vgpt(-)$ by $\H\text{\eighteurm om}\>_{\bold t}\<(\Cal M,-)\set
\drlm\H\text{\eighteurm om}\>\>(\Cal M/\bt\<^n\<\<\Cal M,-)$ \dots}
We work throughout with {\it unbounded\/} complexes, which sometimes
introduces technical complications, but which will ultimately be quite
beneficial in situations involving combinations of right- and
left-derived functors.

Rather than explain further, we simply suggest a perusal of the salient
results---Lemma ~(3.1.1) (especially $(1)\Leftrightarrow(2)$),
(3.1.3)--(3.1.8),
(3.2.3)--(3.2.7). For completeness we have included several results which are
not used elsewhere in this paper. Some readers may prefer going directly
to ~\S4, referring back to~\S3 as needed.
\medskip

\definition{Definition (3.0.1)}
Let $X$ be a topological space and $\O$ a sheaf of commutative rings on~$X$.
A sequence $\bt\set\nmb(t_1,t_2,\dots,t_\mu)$ in $\Gamma(X,\O)$
is {\it proregular\/} if for each $i=\nmb1,2,\nmb\dots,\mu$
and each $r> 0$ there exists an $s>r$ such that in $\O$,
$$
(t_1^s,\dots,t_{i-1}^s)\O: t_i^s\>\subset\>
(t_1^r,\dots,t_{i-1}^r)\O: t_i^{s-r}.
$$

A closed subspace $Z\subset X$ is {\it proregularly embedded in\/ $X$}
if there exists an open covering~$(X_\alpha)_{\alpha\in A}$ of~$X$
and for each $\alpha$
a proregular sequence\vadjust{\kern.6pt}~$\bt_\alpha$ in $\Gamma(X_\alpha,
\O_\alpha)$
(where $\O_\alpha\set\nmb\O|_{X_\alpha})$ such that
$Z\cap X_\alpha$ is the support of $\O_\alpha/\bt_\alpha\O_\alpha\>$.
\enddefinition

\example{Examples} (a) Suppose that $X$ is {\it quasi-compact\/} (not
necessarily Hausdorff,
but every open cover has a finite subcover), and that the $\O$-module~$\O$
is {\it coherent.} Then $\bt$~is proregular if (and clearly only if) for each
$i,r$ as above and each $x\in X\<$,
there exists an $s=s(x)>r$ such that in the stalk $\O_{\<x}\>$,
$$
(t_1^s,\dots,t_{i-1}^s)\O_{\<x\<}: t_i^s\>\subset\>
(t_1^r,\dots,t_{i-1}^r)\O_{\<x\<}: t_i^{s-r}.\tag 3.0.2
$$
Indeed, the ideal sheaves appearing in~(3.0.1) are all coherent, and so
we can take $s(y)=s(x)$ for all points $y$ in some
neighborhood~$W_x$ of ~$x$. If (3.0.2) holds for $s$ then it holds for
all $s'>s$; and since $X$ can be covered\vadjust{\kern.5pt}
by finitely many of the~$W_x\>$, the condition in (3.0.1) is satisfied.

Note that (3.0.2) holds whenever the ring $\O_{\<x}$ is {\it noetherian,}
since then
$$
  (t_1^r,\dots,t_{i-1}^r)\O_{\<x\<}: t_i^{s-r}\jog=\jog
(t_1^r,\dots,t_{i-1}^r)\O_{\<x\<}: t_i^s\quad\text{\,for \,}s\gg r.
$$
Thus if $X$ is quasi-compact, $\O$ is coherent, and all the stalks
$\O_{\<x}$ are noetherian, then every sequence~$\bt$ is proregular.

(b) If (3.0.2) holds, then it also holds when $\O_{\<x}$ is replaced by
any {\it flat\/} $\O_{\<x}$-algebra.
It follows, for example, that if $R$ is a
ring of fractions of a polynomial ring (with any number of indeterminates)
over a noetherian ring, then
every sequence~$\bt$ in $R=\nmb\Gamma(\text{Spec}(\<R), \O_{\text{Spec}(\<R)})$
is proregular; and every closed subscheme~$Z\subset\text{Spec}(\<R)$ such
that $\text{Spec}(\<R)\setminus Z$ is quasi-compact is proregularly embedded.

(c) For an example by Verdier of a non-proregular sequence, and the resulting
homological pathologies, see \cite{I, pp.\,195--198}.
\endexample
\goodbreak

\nextpart{3.1}
Let $(X,\O)$ be as in Definition (3.0.1).
Denote the category of $\O$-modules by $\A$, and let $\D$ be the derived
category of~$\A$. Fix a sequence
$\bt=(t_1,\dots,t_\mu)$ in~$\Gamma(X,\O)$, and set
$$
\bt\<^n\set(t_1^n,\dots,t_\mu^n)\qquad(n>0).
$$
Define the functor $\vgp{\>\bt}\:\A\to\A\>$ by
$$
\vgpt(\G\>)\set\underset{{}^{n>0}}\to
 \drlm\,\sHom_{\O}(\O/\bt\<^n\O,\>\jog \G\>)
\qquad(\G\in\A).
$$
\vskip-1.2\jot
\noindent The stalk of ~$\vgpt(\G\>)$ at any point $x\in X$ is
$$
\vgpt(\G\>)_x=\underset{{}^{n>0}}\to
 \drlm\,\Hom_{\>\O_{\<x}}(\O_{\<x}/\bt\<^n\O_{\<x}\>,\>\jog \G_x)\qquad(x\in
X).
$$
\vskip-1.2\jot
\noindent The (homological) derived functors of $\vgpt$ are
$$
H^i\R\vgpt(\G)=\underset{{}^{n>0}}\to
 \drlm\,\E\eurm x\eurm t_{\O}^i(\O/\bt\<^n\O,\>\jog \G\>)
\qquad(i\ge 0,\ \G\in\A).
$$
\vskip-1.2\jot
\noindent If $\bold s$ is another finite sequence in~$\Gamma(X,\O)$ such that
$\root\of{\bold s\O}=\root\of{\bt\O}$ then
$\vgp{\>\bold s}=\vgp{\>\bt}$.
If~$(X,\O)$ is a scheme and $Z\set \text{Supp}(\O/\bt\O)$ then
$\vgp{\>\bt}=\vgp {Z\>}$, see~(0.1).

\smallskip

For $t\in \Gamma(X,\O)$, let $\K(t)$ be the complex
\smash{$\,\cdots\to 0\to \O\overset
{\vbox to 0pt {\vskip-.85ex\hbox{$\ssize t$}\vss}}\to\to \O\to 0 \to\cdots$}
which in degrees 0 and~1
is multiplication by $t\/\>$ from $\O=:\Cal K^0(t)$ to $\O=:\Cal K^1(t)$,
and which vanishes elsewhere. For $0\le r\le s$, there is
a map of complexes $\K(t^r)\to\nmb \K(t^s)$
which is the identity in degree~0\vadjust{\kern.75pt} and multiplication by~
$t^{s-r}$ in degree~1; and thus we get a direct system of complexes, whose
\smash{$\drlm$}\vadjust{\kern.5pt} we denote by $\KK(t)$.
The stalk\vadjust{\kern.5pt} of~$\KK(t)$ at~$x\in\nmb X$
looks in degrees 0 and~1 like\vadjust{\kern.75pt} the localization map
$\O_{\<x}\to\nmb(\O_{\<x})_t=\nmb\O_{\<x}[1/t]$.
With $\otimes=\smash{\otimes_\O}@,$, set\vadjust{\kern-2\jot}
$$
\spreadlines{1.5\jot}
\gather
\K(\bt)\set \K(t_1)\otimes\dots
 \otimes \K(t_\mu),\\
\KK(\bt)
 \set \underset {{}^{n>0}}\to \drlm\> \K(\bt\<^n)
  = \KK(t_1)\otimes\dots \otimes \KK(t_\mu);
\endgather
$$
and for any complex~$\F$ of $\O$-modules set
$$
\K(\bt,\F\>)\set \K(\bt)\<\otimes\<\F\>,\qquad
 \KK(\bt,\F\>)\set \KK(\bt)\<\otimes\<\F.
$$
Since the complex~$\KK(\bt)$ is flat and bounded, the functor
of complexes
$\KK(\bt,-)$ takes quasi-isomorphisms to quasi-isomorphisms
\cite{H, p.\,93, Lemma 4.1, b2},
and so may be regarded as a functor from ~$\D$ to $\D$.

After choosing a quasi-isomorphism~$\varphi$ from $\F$ to a
K-injective $\O$-complex~$\Cal L^\bullet$ \cite{Sp,~p.\,138, Thm.\,4.5},
we can use the~natural\- identifications
$$
\vgpt(\Cal L^j)=\drlm \ker\bigl(\Cal K^0(\bt\<^n\!,\>\Cal L^j)\to
  \Cal K^1(\bt\<^n\!,\>\Cal L^j)\bigr)
 =\ker\bigl(\Cal K_{\<\<\sssize\infty}^0\<\<(\bt,\Cal L^j)\to
 \Cal K_{\<\<\sssize\infty}^1\<\<(\bt,\Cal L^j)\bigr)
\quad(j\<\in\<\Bbb Z)
$$
to get a $\D$-morphism
$$
\nopagebreak
\delta'=\delta'(\F\>)\:\R\vgpt(\F\>)\cong\vgpt(\Cal L^\bullet)\hookrightarrow
 \KK(\bt,\Cal L^\bullet)\cong\KK(\bt,\F\>),
$$
easily checked to be functorial in $\F$ (and in particular, independent
of~$\varphi$).

In proving the next Lemma, we will see that proregularity of $\bt$ implies that
$\delta'(\F\>)$ is always an isomorphism. And the converse holds if cohomology
on $X$ commutes with filtered direct limits,
for example if $X$ is {\it compact\/} (i.e.,
quasi-compact and Hausdorff) \cite{Go, p.\,194, Thm.\,4.12.1}, or if $X$
is {\it quasi-noetherian\/} \cite{Ke, p.\, 641, Thm.\,8}.
Kempf defines $X$ to be quasi-noetherian
if its topology has a base of quasi-compact open sets,
if the intersection of any
two quasi-compact open subsets of ~$X$ is again quasi-compact, and if
$X$ itself is quasi-compact. We prefer to use the term  {\it concentrated\/}.
For example, if $X$ is noetherian (i.e., every open subset is quasi-compact)
then $X$ is concentrated.
A scheme is concentrated iff it is quasi-compact and quasi-separated
\cite{GrD, p.\,296, Prop.\,(6.1.12)}.
\footnote
{where, for the implication $\text{d)}\Rightarrow \text{a)}$, the family
$(U_{\<\alpha})$ should be a base of the topology.%
}

\proclaim{Lemma (3.1.1)} Let\/ $\bt=(t_1,\dots,t_\mu)\ \;
 \bigl(t_i\in \Gamma(X,\O)\bigr)$ and $\delta'$ be as above, and
suppose~that\/ $X$ is compact or concentrated. Then
the following are equivalent:
\smallskip
$(1)$ The sequence\/ $\bt$ is proregular\/
$($Definition\/ $(3.0.1)).$

$(2)$ \kern1pt For any\/~$\F\in\D,$\ the map\/
$\delta'(\F\>)\:\R\vgpt(\F\>)\to \KK(\bt,\F\>)$
is an isomorphism.\vadjust{\kern.75pt}

$(2)'$ \kern-1.2pt
For any injective\/ $\O$-module\/ $\J$ and every\/ $i\ne0,$\
$H^i\KK(\bt,\J\>)=0$.

$(3)$ \kern1.5pt For any flat\/ $\O$-module\/ $\P$ and every\/ $i\ne 0,$\ the
inverse system\/
$$
\bigl(H_i(\bt\<^r\<, \P)\bigr)_{\<r>0}\>\>\set
\bigl(H^{-i}\>\sHom_\O(\K(\bt\<^r), \P)\bigr)_{\<r>0}
$$
is essentially null, i.e., for each\/~$r$ there is an\/
$s>r$ such that the natural map\/
$H_i(\bt\<^s\<, \P)\to H_i(\bt\<^r\<,\P)$
is the zero map\vadjust{\kern1pt}.

$(3)'$ \kern1pt For every\/ $i\ne 0$ the inverse system\/
$\bigl(H_i(\bt\<^r\<,\O)\bigr)_{\<r>0}$
is essentially null.

$(3)'{}'$ \kern-1.5pt
The inverse system\/ $\bigl(H_1(\bt\<^r\<, \O)\bigr)_{r>0\>}$ is
essentially null.\vadjust{\kern1pt}

$(4)$ \kern3.3pt Setting for any\/ $\O$-complex\/ $\F\<,$\ and\/ $r>0,$
$$
\H_{@,@,\bt\<^r }\>\F\set\sHom(\O/\bt^r\O,\jog\F\>)
 =\ker\bigl(\Cal K^0\<(\bt^r\<,\jog\F\>)\to
\Cal K^1\<(\bt^r\<,\jog\F\>)\bigr),
$$
we have that if\/ $\Cal R$ is a complex of\/ $\O$-injectives and\/
$\Cal C_r$ is the mapping cone of\/
$
\H_{@,@,\bt\<^r}\<\Cal R \hookrightarrow\nmb \K(\bt\<^r\<\<,\Cal R)
$
then for each\/~$r>0$ there is an\/ $s>r$ such that the natural map\/
$\>\Cal C_r\to\Cal C_{s}$ vanishes in\/~$\D$.
\endproclaim
\goodbreak
\proof
We proceed as follows. \vadjust{\kern.5\jot}

(A):\quad $(1) \Rightarrow  (3)' \Rightarrow (3)'{}' \Rightarrow (1)$.

(B):\quad $(3) \Rightarrow (3)' \Rightarrow (2)' \Leftrightarrow (2)
 \Rightarrow (3)$.

(C):\quad $(1) \Rightarrow (4) \Rightarrow (2)'$.\vadjust{\kern2pt}

\noindent(Substituting $\P$ for~$\O$ in the proof of  $(1)\Rightarrow (3)'$
yields a direct proof of~$(1)\Rightarrow\nmb (3)$.)

\noindent The hypothesis ``$X$ compact or concentrated"
will be needed only for $(2)\Rightarrow(3)$.
\medskip
{\bf(A).}
Assuming (1), we prove $(3)'$ by induction on $\mu$. For $\mu =1$,
the assertion amounts to the vanishing (in~$\O$) of
$\>t_1^{s-r}(0\quot t_1^s)$ when $s\gg r$,
which we get by taking $i=1$ in Definition~(3.0.1). For $\mu>1$,
there is an obvious direct system of split exact sequences of complexes
$$
0\to \O_{\<r}'[-1]\to \Cal K^\bullet(t_\mu^r) \to \O_r\to 0\qquad(r>0)
\tag 3.1.1.1
$$
where $\O_{\<r}'\>\set \O =: \O_r$ for all~$r$, where the map $\O_{\<s}'\to
\O_r'\ (s>r)$
is multiplication by~$t_\mu^{s-r}\<$, and $\O_{\<s}\to \O_r$ is the identity.
{}From this system we derive an inverse system of exact sequences
$$
\spreadlines{1\jot}\multlinegap{2pt}
\multline
0\to\sHom_\O(\Cal K^\bullet((t_1^r,\dots,t_{\mu-1}^r))\>\>\otimes
\>\>\O_r\>,\>\O)\to
 \sHom_\O(\Cal K^\bullet((t_1^r,\dots,t_{\mu-1}^r))\>\>\otimes\>\> \Cal
K^\bullet(t_\mu^r),\>\O)\\
  \to\sHom_\O(\Cal K^\bullet((t_1^r,\dots,t_{\mu-1}^r))
    \<\otimes \<\O_{\<r}'\>,\>\O)[1]\to 0
\endmultline
$$
whence an inverse system of exact homology sequences, with
$\,\Cal I^{\jog[r]}\set(t_1^r,\dots,t_{\mu-1}^r)\O$,
$$
\def\pt{\phantom{t_\mu^r}}
\spreadlines{1\jot}
\multline
\cdots@>\pt>> H_i\bigl((t_1^r,\dots,t_{\mu-1}^r), \O\bigr)
 @>t_\mu^r>>H_i\bigl((t_1^r,\dots,t_{\mu-1}^r), \O\bigr)
  @>\pt>> H_i\bigl(\bt\<^r,\O\bigr)\\
\ \quad\qquad@>\pt>> H_{i-1}\bigl((t_1^r,\dots,t_{\mu-1}^r), \O\bigr)
  @>t_\mu^r>>H_{i-1}\bigl((t_1^r,\dots,t_{\mu-1}^r), \O\bigr)
    @>\pt>>\cdots\hskip30pt  \\
\cdots@>\pt>> H_1\bigl((t_1^r,\dots,t_{\mu-1}^r), \O\bigr)
  @>\pt>> H_1\bigl(\bt\<^r,\O\bigr)
   @>\pt>> (\Cal I^{\jog[r]}\quot t_\mu^r)/\Cal I^{\jog[r]}
    @>\pt>> 0.
\endmultline
$$
Now the inductive hypothesis quickly reduces the problem
to showing that the inverse system
$T_r\set(\Cal I^{\jog[r]}\quot t_\mu^r)/\Cal I^{\jog[r]}$,
with maps $T_s\to T_r\ (s>r)$ given by multiplication by~$t_\mu^{s-r}\<$, is
essentially null;
and that results from Definition (3.0.1) with~$i=\nmb\mu$.

Thus (1) implies $(3)'$, of which $(3)'{}'$ is a special case.
Conversely, assuming $(3)'{}'$ we~get~(1) from the surjections (as above)
$$
H_1\bigl((t_1^r,\dots,t_i^r),\O\bigr)
   \twoheadrightarrow \bigl((t_1^r,\dots,t_{i-1}^r)\O: t_i^r\bigr)/
     (t_1^r,\dots,t_{i-1}^r)\O   \qquad (1\le i\le\mu).
$$

{\bf(B).} Take $\P=\O$ in (3) to get $(3)'$.
If $\J$ is an injective $\O$-module, then
$$
\spreadlines{1\jot}
\align
H^i\KK(\bt,\jog\J\>)=
H^i\>\> \underset {{}^{r>0}}\to \drlm \K(\bt\<^r\<,\jog\J\>))
&\cong\underset{{}^{r>0}}\to
   \drlm H^i\>\sHom\bigl(\sHom(\K(\bt\<^r),\O),\jog\J\>\bigr)\\
&=\underset {{}^{r>0}}\to
   \drlm \sHom\bigl((H_{-i}(\bt\<^r\<,\O)),\jog\J\>\bigr)
\endalign
$$
and consequently $(3)'\Rightarrow(2)'$.

$(2)'$ implies, for any $\O$-complex~$\F$,
that if $\F\to\Cal L^\bullet$ is a quasi-isomorphism
with $\Cal L^\bullet$ both K-injective and injective \cite{Sp, p.\,138, 4.5},
then
the $j$-th column $\KK(\bt,\Cal L^j)$ of the double complex
$(\Cal K_{\<\<\sssize\infty}^i(\bt)\otimes\nmb
  \Cal L^j)_{0\le i\le \mu,\;j\in\>\lower.1ex\hbox{$\ssize\Bbb Z$}}$
is a finite resolution of~$\vgpt(\Cal L^j)$, so that\vadjust{\kern1pt}
the inclusion
$\vgpt(\Cal L^\bullet)\hookrightarrow\KK(\bt,\Cal L^\bullet)$ is a
quasi-isomorphism;  and (2)~follows.
Conversely, since $\R\vgpt(\J)\cong \vgpt(\J)$, $(2)'$ follows from~(2).
\goodbreak
To deduce (3) from~$(2)'$ we imitate \cite{Gr, p.\,24}.
There exists a monomorphism of the $\O$-module $H^i\sHomb(\K(\bt\<^r), \P)$
into an injective $\O$-module~$\J'\<$, giving rise naturally to an element of
$$
\underset{^{s>r}}\to\drlm\>
 \Hom\bigl(H^i\sHomb(\K(\bt\<^s), \P),\jog\J'\bigr);
 \tag 3.1.1.2
$$
and it will suffice to show that this element is zero. Noting that
homology commutes with the exact functor~$\Hom(-,\J')$ and
with~$\drlm$\kern-.75pt,
noting that $\K(\bt\<^s)$ is a finite-rank free $\O$-complex,
setting $\Gamma(-)\set\Gamma(X,-)$,
and setting $\J\set\nmb\sHomb(\P,\jog\J')$ (which is
an injective $\O$-module since $\P$ is flat),
we can rewrite (3.1.1.2) as
$$
\align
&\quad\ H^i\drlm\Gamma\>\sHomb\bigl(\sHomb(\K(\bt\<^s), \P),\jog\J' \bigr) \\
&=H^i\drlm\Gamma\>\sHomb\bigl(\sHomb(\K(\bt\<^s), \O)\otimes\P,\jog\J'\bigr)\\
&=H^i\drlm\Gamma\>\sHomb\bigl(\sHomb(\K(\bt\<^s), \O),\jog\J \bigr)\\
&=H^i\drlm\Gamma\bigl(\K(\bt\<^s)\otimes\J \bigr),
\endalign
$$
or again,
since $\Gamma$ commutes with \smash{$\drlm$}\vadjust{\kern1.3pt} ($X$ being
compact or concentrated),
as $H^i@,@,\Gamma@,@,\KK(\bt,\jog\J)$.
But by~$(2)'\<$, $\KK(\bt,\jog\J)$ is  a resolution of
$\vgp{Z}\J\<$, and as a \smash{$\drlm$}\vadjust{\kern1.3pt}
of injective complexes, is a
complex of $\Gamma$-acyclic sheaves (since $H^i(X,-)$
commutes\vadjust{\kern1.3pt}
with~\smash{$\drlm$}); also $\vgp{Z}\J\<$, the \smash{$\drlm$} of
the flabby sheaves $\sHom(\O/\bt\<^n\O,\jog\J)$, is $\Gamma$-acyclic;
and so
$$
H^i@,@,\Gamma@,@,\KK(\bt,\jog\J)=H^i(X,\vgp{Z}\J)=0\qquad (i\ne 0).
$$

\smallskip
{\bf(C).}
$(2)'$ follows from (4) (with $\Cal R\set\J\>$) upon application of
$\smash{\drlm}$ to the direct system $(E_s)_{s>0}$ of exact sequences
$$
\cdots\to H^{i-1}\bigl(\Cal C_s  \bigr)
 \to  H^i\bigl(\H_{@,@,\bt ^s}\<\Cal J\>\bigr)
 \to H^i\bigl( \K(\bt\<^s\<,\Cal J\>) \bigr)
 \to H^i\bigl(\Cal C_s  \bigr)\to\cdots
\tag "$(E_s)$"
$$

Finally, $(1)\Rightarrow(4)$ is proved by induction on~$\mu$.
Consider the following
natural commutative diagram of complexes, in~which $s>r$ and $\nu\set\mu-1$.
$$
\eightpoint
\minCDarrowwidth{.18in}
\def\1{\H_{(t_1^r,\dots,t_{\nu}^r)}\H_{(t_{\<\mu}^r)}\Cal R}
\def\2{\H_{(t_1^r,\dots,t_{\nu}^r)}\H_{(t_{\<\mu}^s)}\Cal R}
\def\3{\H_{(t_1^s,\dots,t_{\nu}^s)}\H_{(t_{\<\mu}^s)}\Cal R}
\def\4{\H_{(t_1^r,\dots,t_{\nu}^r)}\K((t_{\<\mu}^r),\Cal R)}
\def\5{\H_{(t_1^r,\dots,t_{\nu}^r)}\K((t_{\<\mu}^s),\Cal R)}
\def\6{\H_{(t_1^s,\dots,t_{\nu}^s)}\K((t_{\<\mu}^s),\Cal R)}
\def\7{\K\bigl((t_1^r,\dots,t_{\nu}^r),\K((t_{\<\mu}^r),\Cal R)\bigr)}\<\<
\def\8{\K\bigl((t_1^r,\dots,t_{\nu}^r),\K((t_{\<\mu}^s),\Cal R)\bigr)}\<\<
\def\9{\K\bigl((t_1^s,\dots,t_{\nu}^s),\K((t_{\<\mu}^s),\Cal R)\bigr)}
\CD
\1@>>>\2@>>>\3 \\
@VaVV @VVbV @VVcV  \\
\4@>>>\5@>>>\6 \\
@VdVV @VVeV @VVfV  \\
\7@>>>\8@>>>\9
\endCD
$$
\vskip1pt\noindent
(When $\mu=1$, interpret $d$, $e$, and $f$ to be the
identity map of $\K\bigl((t_1^{\sssize\bullet}),\Cal R\bigr)$.)
We need to prove, modulo obvious identifications, that the natural
$\D$-homomorphism
$\Cal C_{da}\to\Cal C_{\<f\<c}$ from the cone of~$da$
to the cone of~$f\<c$ is zero if $s\gg r$. That's true---obviously
when $\mu=1$ and by the inductive
hypothesis otherwise---for $\>\Cal C_e\to \Cal C_{\<f}\>$, hence for
$\>\Cal C_d\to \Cal C_{\<f}\>$.
It's also true for $\>\Cal C_a\to \Cal C_c\>$:
indeed, $a$ being a monomorphism, $\Cal C_{a[1]}$ is $\D$-isomorphic to
$\operatorname{coker}(a[1])$,
which is isomorphic as an $\O$-complex to the mapping cone
of the inclusion $\iota\:t_\mu^r\H_{(t_1^r,\dots,t_{\nu}^r)}\Cal R
\hookrightarrow \H_{(t_1^r,\dots,t_{\nu}^r)}\Cal R$, which cone is
itself $\D$-isomorphic to $\operatorname{coker}(\iota$); working
through these isomorphisms we can identify the
$\D$-map $\Cal C_a\to \Cal C_c$ with
$$
\smash{\bigl(\H_{(t_1^r,\dots,t_{\nu}^r)}\Cal R/
  t_{\<\mu}^r\H_{(t_1^r,\dots,t_{\nu}^r)}\Cal R\bigr)[-1]
@>\pm t_{\<\mu}^{s-r}>>
\bigl(\H_{(t_1^s,\dots,t_{\nu}^s)}\Cal R/
  t_{\<\mu}^s\H_{(t_1^s,\dots,t_{\nu}^s)}\Cal R\bigr)[-1];}
$$
and recalling Definition (3.0.1), we conclude then  via the
next Lemma---with $t=t_{\<\mu}$, $\I=\nmb(t_1^s,\dots,t_{\nu}^s)\O$,
and $\I'=(t_1^r,\dots,t_{\nu}^r)\O$.

\proclaim{Lemma (3.1.2)}
Let $\I,$\ $\I'$ be $\O$-ideals, let $t\in\Gamma(X, \O),$\ and let $r\le s$ be
integers
such that $t^{s-r}(\I\<\<\quot t^s)\subset \I'\<\<$. Then for any injective\/
$\O$-module~$\J$ and any open\/ $U\subset X$
we have, setting\/ $\G_U\set\G|_U$\ for any\/ $\O$-module~$\G\>$:
$$
t^{s-r}\Hom(\O_U/\I_U',\> \J_U)\subset
t^s\Hom(\O_U/\I_U,\> \J_U).
$$
\endproclaim

\proof
For any map $\alpha\:\O_U/\I_U'\to \J_U$,
the kernel of $\smash{\O_U@>t^{s}>> \O_U/\I_U}$
annihilates $t^{s-r}\alpha$
(because $(\I\<\<\quot t^s)t^{s-r}\subset \I'$),
and so there is an $\O_U$-homomorphism
$$
\psi=\psi_{r\<,s\<,\alpha}\:t^s(\O_U/\I_U\>)\to
 \sHom(\O_U/\I_U',\jog \J_U)\subset\J_U
$$
with $\psi(t^s+\I_U\>)=t^{s-r}\alpha$.  Since $\J_U$ is
an injective $\O_U$-module, $\psi$ extends to a map~
$\psi^0\:\O_U/\I_U\to\nmb\J_U\>$, and then
$$
t^{s-r}\alpha=\psi^0(t^s+\I_U\>)=t^s\psi^0(1+\I_U\>)
\in t^s\Hom(\O_U/\I_U,\jog \J_U)\subset\Gamma(U,\>\J_U).
$$
\vskip-2\jot
\rightline{\hfill\qed}
\enddemo
\goodbreak

Now there is a natural commutative diagram, with $\Hom\set\Hom_{\D}\>$,
$$
\CD
\Hom(\Cal C_{da},\jog \Cal C_a)@>>>\Hom(\Cal C_{da},\jog \Cal
C_{da})@>>>\Hom(\Cal C_{da},\jog \Cal C_d)\\
@VVV @VV\rho V @VVV \\
\Hom(\Cal C_{da},\jog \Cal C_c)@>\sigma>>\Hom(\Cal C_{da},\jog \Cal
C_{\<fc})@>>>\Hom(\Cal C_{da},\jog \Cal C_{\<f})
\endCD
$$
in which the rows are {\it exact\/} (by \cite{H, p\,23, Prop.\,1.1\,b)} and
the octahedral axiom [{\it ibid.,}~p.\,21]).
For fixed
$r$ and variable~$s$, the outside columns
form---as we have just seen---essentially null direct systems,
whence so does the middle column. The desired conclusion results.

\smallskip
This completes the proof of Lemma (3.1.1).
\endproof
\smallskip
\smallskip
With no assumption on the topological space $X$ we define as in~(2.1) {\it
mutatis mutandis} a functorial map
$$
\botsmash{\psi_\bt'(\E,\F\>)\:\E\Otimes\R\vgpt(\F\>)\iso
\R\vgpt(\E\Otimes\F\>)}
\qquad(\E,\F\in\D).
$$

\proclaim{Corollary (3.1.3)}
 \!If\/ $\bt$ is proregular then\/ $\>\psi_\bt'(\E,\F\>)$ is an isomorphism for
all\/ $\E,\>\F$.
\endproclaim
\proof
Assume, as one may, that $\E$ is K-flat, and check that
the following diagram---whose bottom row is the natural
isomorphism---commutes:
$$
\CD
\E\otimes\R\vgpt(\F\>) @>\psi_\bt'(\E,\>\F\>)>>\R \vgpt(\E\otimes\F\>) \\
@V \text{via }\delta'(\F\>) VV @VV\delta'(\E\otimes\F\>)V  \\
\E\otimes\KK(\bt,\F\>)
@>\Iso>\phantom{\psi_\bt'(\E,\F\>)}>\KK(\bt,\E\otimes\F\>)
 \endCD
$$
\vskip-6pt\noindent
By the implication $(1)\Rightarrow(2)$ in~Lemma~(3.1.1)
(whose proof did not need $X$ to be compact or concentrated), the maps
$\delta'(\F\>)$ and $\delta'(\E\otimes\F\>)$ are also isomorphisms,
and the assertion follows.
\endproof

\proclaim{Corollary (3.1.4)}
If\/ $\bt$ and\/ $\bt^*$ in\/~$\Gamma(X,\O)$  are such
that\/ $\bt^*$ and\/ $(\bt,\bt^*)$ are both proregular, then the natural map\/
$\R\vgp{(\bt,\bt^*\<)}\to\R\vgp\bt\smcirc\R\vgp {\bt^*}$ is an isomorphism.
\endproclaim

\proof
Proregularity of $(\bt,\bt^*)$ trivially implies that of~$\bt$ (and also,
when $X$ is compact or concentrated, of~$\bt^*$, see remark preceding (3.1.6)
below). By (3.1.1)(2), the assertion results from the equality
$\KK((\bt,\bt^*),-)=\KK(\bt,\jog\KK(\bt^*,-))$.
\endproof

\proclaim{Corollary (3.1.5)}
Let\/ $(X,\O)$ be a scheme and\/ $Z\subset X$ a
proregularly embedded subscheme.

{\rm(i)} The map
$
\botsmash{\psi'\:\E\Otimes\R\vgp Z(\F\>)\iso \R\vgp Z(\E\Otimes\F\>)}
$
of\/ $(2.1)$ is an isomorphism for all $\E,\F \in\D(X)$.\vadjust{\kern.6pt}

{\rm(ii)} If\/ $Z^*\subset X$ is a closed subscheme such that\/ $Z^*$ and\/
$Z\cap Z^*$ are both proregularly embedded, then the natural functorial map\/
$\R\vgp{Z\cap Z^*}\to\R\vgp Z\smcirc\R\vgp {Z^*}$ is an isomorphism.

{\rm(iii)} $\R\vgp Z\bigl(\D\qcd(X)\bigr)\subset \D\qcd(X)$.
\endproclaim

\proof
The assertions are essentially local on~$X\<$, so the first two follow from
(3.1.3) and~(3.1.4) respectively, and the third from (3.1.1)(2), see
\cite{H, p.\,98, Prop.\,4.3}.
\endproof

Assume now that $X$ is compact or concentrated.
If $\bt\<^*$ is a permutation of~$\bt$ then there is an obvious functorial
isomorphism $\KK(\bt\<^*,-)\iso \KK(\bt,-)$, and so by~
Lemma~(3.1.1)(2), $\bt\<^*$ is proregular $\Leftrightarrow$ so~is~$\bt$.
More generally:

\proclaim{Corollary (3.1.6)}
Let\/ $\bt=(t_1,\dots,t_\mu)$ be, as before, a sequence in\/ $\Gamma(X,\O),$\
with\/ $X$ compact or concentrated,
and let\/ $\bt\<^*\set(t_1^*,\dots,t_\nu^*)$ be a sequence\vadjust{\kern1pt}
in\/
$\Gamma(X,\,\root\of{\bt\O}\>\>)$.
Then the sequence\/
$(\bt\<^*\<,\bt)\set(t_1^*,\dots,t_\nu^*\>,t_1,\dots,t_\mu)$
is proregular $\Leftrightarrow$ so is\/ $\bt$.
In~particular, if\/~$\>\>\root\of{\bt\<^*\O}=\nmb\root\of{\bt\O}$ then\/
$\bt\<^*$ is proregular~$\Leftrightarrow$ so is\/~$\bt$.
\endproclaim

\proof
It suffices to treat the case $\nu=1$. Since (clearly)
$\vgp{(\bt\<^*\!,\jog\bt)}=\vgpt\>$, and in view of~(3.1.1)(2), we need
only show, with $t\set t_1^*\>$, that for any $\O$-complex~$\F\>$ the natural
functorial map
$$
\KK((t,t_1,\dots,t_\mu),\F\>)=\KK(t)\otimes\KK(\bt,\F\>)
\to \O\otimes\KK(\bt,\F\>)=\KK(\bt,\F\>)
$$
induces homology isomorphisms. The kernel of this degreewise split
surjective map is $\O_t[-1]\otimes\KK(\bt,\F\>)$, where $\O_t$ is the
direct limit of the system $(\O_n)_{n>0}$ with $\O_n\set\nmb \O\>$
for all~$n$ and with $\O_r\to \O_s\ (r\le s)$  multiplication by~$t^{s-r}$; and
it~will suffice to show
that this kernel is {\it exact,} i.e., that {\it for $j\in\Bbb Z$
and\/ $r>0,$\ any~section of\/~
$H^j(\KK(\bt\<^r,\F\>))$ over an open\/ $U\subset X$
is locally\vadjust{\kern.6pt} annihilated by a power of\/~$t$.}
Since~$t\in\nmb\root\of{\bt\O}$ we can
replace~$t\/$ by~$t_i\ (1\le i\le\mu)$ in this last statement, whereupon
it becomes well-known---and easily proved
by induction on~$\mu$, via~(3.1.1.1).
\endproof

\proclaim{Corollary (3.1.7)}
Let\/ $(X,\O)$ be a quasi-separated scheme and\/
$Z\subset X$ a proregularly embedded subscheme.
If\/ $X_0\subset X$ is a quasi-compact open subset,
$\O_0\set\nmb \O|_{X_0\>},$\ and\/ $\bt_0$~is a finite sequence
in\/ $\Gamma(X_0, \O_0)$ such that\/
$Z\cap X_0$ is the support of\/~$\O_0/\bt_0\O_0\>,$\ then
$\bt_0$~is~proregular.
\endproclaim

\proof
$X_0$ is covered by finitely many of the open sets~$X_0\cap X_\alpha$ with
$X_\alpha$ as in Definition (3.0.1), and we may assume that each $X_\alpha$ is
quasi-compact, whence so is $X_0\cap X_\alpha$ (since $X$ is quasi-separated).
So it suffices to apply (3.1.6) to~
$X_0\cap\nmb X_\alpha\>$, with $\bt\set \bt_0$ and $\bt\<^*\set \bt_\alpha$.
\endproof

\eightpoint
Let $(X,\O)$ be a scheme, let $j\:\A\qcd=\A\qcd(X)\hookrightarrow\A$ be the
inclusion of the category of quasi-coherent $\O$-modules into the category of
all $\O$-modules, and let
$\boldkey j\:\D(\A\qcd)\to\D(\A)=:\D$ be
the~corresponding derived-category functor.

\proclaim{Proposition (3.1.8)}
If\/ $(X,\OX\<)$ is a quasi-compact separated scheme and\/ $Z\subset X$ is
proregularly embedded, then the functor
$$
\varGamma\>\qcu_{\mkern-8mu Z}\set
 \vg Z\smcirc j=\vgp Z\smcirc j\:\A\qcd\to\A\qcd
$$
has a derived functor
$$
\R\varGamma\>\qcu_{\mkern-8mu Z}\:\bold  D(\A\qcd)\to\D(\A\qcd);
$$
and the natural functorial map\/
$\>\boldkey j\smcirc\R\varGamma\>\qcu_{\mkern-8mu Z}
 \to \R\vgp Z\smcirc \boldkey j$ is an isomorphism.
\endproclaim

\remark{Remark}
For quasi-compact separated $X\<$, $\>\boldkey j\jog$
induces an {\it equivalence of categories\/} from $\D(\A\qcd)$ to
$\D\qcd(X)$ \cite{BN,~p.\,230, Cor.\,5.5} (or see (1.3) above).
Therefore any $\F\in\D\qcd(X)$
is isomorphic to a quasi-coherent complex. In this case, then,
(3.1.8)~embellishes assertion (iii) in~(3.1.5). (The following proof
does not, however, depend on \cite{BN} or (1.3).)
\endremark\smallskip

Proposition (3.1.8) is a consequence of:

\proclaim{Lemma (3.1.8.1)}
For any inclusion\/ $i\:U\hookrightarrow X$ with\/
$U$ affine open, and any\/~$\J$ which is injective in\/~$\A\qcd(U),$\
the natural map\/ $\vgp Z(i_*\<\J\>)\to\nmb\R\vgp Z(i_*\<\J\>)$ is a
$\D$-isomorphism.
\endproclaim

Indeed, if $\>\G\in\A\qcd$, if $(U_\a)_{1\le\a\le n}$ is an affine
open cover of~$X$,
with inclusion maps $i_\a\:U_\a\hookrightarrow X\<$, and if for each~$\a$,
$i_\a^*\G\to\J_\a$ is a monomorphism with $\J_\a$ injective in~$\A\qcd(U_\a)$,
then $i_{\a*}\>\J_\a$ is $\A\qcd$-injective
(since $i_{\a*}\:\A\qcd(U_\a)\to\A\qcd$ has an exact left adjoint),
and there are obvious monomorphisms
$\G\to\oplus_{\a=1}^n \>i_{\a*}\>i_\a^*\G\to \oplus_{\a=1}^n\> i_{\a*}\>\J_\a.$
Thus the category $\A\qcd$ has enough injectives; and since, by (3.1.8.1),
$$
\varGamma\>\qcu_{\mkern-8mu Z}(\oplus_{\a=1}^n\> i_{\a*}\>\J_\a)\smcong
\oplus_{\a=1}^n\>\R\vgp Z(i_{\a*}\>\J_\a),
$$
and the functor $\R\vgp Z$ is bounded above and below (by Lemma~(3.1.1)(2) and
quasi-compactness of~$X\<$), it follows from ~\cite{H, p.\,57, $\gamma$b}
and its proof that $\R\varGamma\>\qcu_{\mkern-8mu Z}$ exists and is
bounded above and below. And then the isomorphism assertion in~(3.1.8)
follows from  \cite{H, p.\,69, (iii) and (iv)}.

It remains then to prove Lemma (3.1.8.1).

Since $X$ is concentrated, there is a finite-type $\OX$-ideal~$\I$
such that $Z=\text{Supp}(\OX/\I\>)$. With $\O_U\set i^*\OX@!$, $\>\I_U\set
i^*\I$,
we have for any $\O_U$-module~$\E$,
$$
\align
\vgp Z @,i_*(\E) &= \drlm\H\text{\eighteurm om}\>(\OX/\I^{\>n}\<,\jog i_*\E) \\
&=\drlm i_*\H\text{\eighteurm om}\> (\O_U/\I_U^n,\jog \E) \\
&=i_*\jog\smash{\drlm}\H\text{\eighteurm om}\>
 (\O_U/\I_U^n,\jog\E)=i_*\vgp{Z\cap U}(\E)
\endalign
$$
where the interchange of $\jog\smash{\drlm}$ and $i_*$ is justified by
\cite{Ke, p.\,641, Prop.\,6}. Since the map~$i$ is affine, and
$i_*$ takes $\O_U$-injectives
to $\OX$-injectives, and since for any $\O_U$-injective ~$\Cal L$,
$\vgp {Z\cap U}(\Cal L)$ is a $\smash{\drlm}$ of flabby sheaves and hence
$i_*$-acyclic \cite{Ke, p.\,641, Cors.\, 5 and~7}, therefore
$$
\R\vgp Z(i_*\<\J)\cong\R\vgp Z(\R\>i_*\<\J)\cong \R(\vgp Z i_*)(\J)
=\R(i_*\vgp{Z\cap U})(\J)=\R\>i_*\R\vgp{Z\cap U}(\J).
$$
Referring again to the ring-theoretic analogue of $(3.1.1)(2)'$
\cite{Gr, p.\,24, Lemme~9,\;b)}, we see that
$\R\vgp{Z\cap U}(\J)\cong\vgp{Z\cap U}(\J)$; and since $i$ is affine
and $\vgp{Z\cap U}(\J)$ is quasi-coherent, therefore
$$
\R\>i_*\R\vgp{Z\cap U}(\J)\cong\R\>i_*\vgp{Z\cap U}(\J)\cong i_*\vgp{Z\cap
U}(\J)\cong \vgp Z i_*(\J),
$$
whence the desired conclusion.\qed
\par\tenpoint

\nextpart{3.2} The map
$$\delta'=\delta'(\F\>)\:\R\vgpt(\F\>)\to\KK(\bt,\F\>)\qquad
 \bigr(\F\in\D\bigl)
$$
remains as in \S3.1. Let $Z$ be the support of
$\O/\bt\O$, a closed subset of~$X$.
In the following steps a)--d), {\it we construct a functorial map}
$$
\delta=\delta(\F\>)\:\KK(\bt,\F\>)\to \R\vg{Z}(\F\>)\qquad
 \bigr(\F\in\D\bigl)
$$
{\it such that\/ $\delta\smcirc\delta'\:\R\vgpt(\F\>)\to \R\vg{Z}(\F\>)$
coincides with the map induced by the obvious inclusion\/}
$\vgpt\hookrightarrow\vg Z\>$.
\smallskip
a) As in the definition of $\delta'$  we may assume that $\F$  is
K-injective, and injective as well
(i.e., each of its component $\O$-modules
$\F^n\ (n\in\Bbb Z)$ is injective) \cite{Sp, p.\,138,~4.5}.  If
$\smash{U\set (X\setminus Z)\overset
  {\vbox to 0pt {\vskip-.75ex\hbox{$\ssize i$}\vss}}
 \to\hookrightarrow X}$
is the inclusion map,\vadjust{\kern-.5pt}
then the canonical sequence of complexes
$0\to\vg Z(\F\>)\hookrightarrow\F@>\eta>> i_*i^*\F\to 0$ is exact,
and\vadjust{\kern.5pt}
{\it there results a natural quasi-isomorphism $\vg Z(\F\>)\to \Cal C_\eta[-1]$
where $\Cal C_\eta$ is the cone of $\eta$.}\vadjust{\kern1.5pt}

b) Let $\Cal K_\flat$ be the complex\vadjust{\kern1pt}
$\Cal K_{\<\<\sssize\infty}^1\<\<(\bt)\to
 \Cal K_{\<\<\sssize\infty}^2\<\<(\bt)\to\dots\;$
($\Cal K_\flat^0\set\Cal K_{\<\<\sssize\infty}^1\<\<(\bt)$,
$\Cal K_\flat^1\set\Cal K_{\<\<\sssize\infty}^2\<\<(\bt),\> \dots)$
There is an obvious map of complexes
$\O\set\Cal K_{\<\<\sssize\infty}^0\<\<(\bt)\to\Cal K_\flat\>$,
inducing\vadjust{\kern1pt} for any
complex~$\F$ a map $\xi=\xi(\F\>)\:\F=\O\otimes\F\to\Cal K_\flat\<\otimes\F$,
{\it whose cone $\Cal C_\xi$ is\/} $\KK(\bt,\F\>)[1]$.\vadjust{\kern1pt}

c) Since $\bt\O_U=\nmb\O_U\ (\O_U\set\O|_U)$, the complex $i^*\KK(\bt)$ is
homotopically
trivial at each point of~$U\<$, and hence for any $\F$ the complex
$i^*\KK(\bt,\F\>)$ is exact. In other words,
$i^*\xi(\F\>)\:i^*\F\to i^*\Cal K_\flat\<\<\otimes\nmb i^*\F$
is a {\it quasi-isomorphism\/} for all~$\F$.

Let $\sigma\:i^*\Cal K_\flat\<\<\otimes i^*\F\to\Cal L\>$ be a
quasi-isomorphism
with $\Cal L\jog$ K-injective. Then $\sigma\smcirc i^*\xi\:i^*\F\to\Cal L$ is a
quasi-isomorphism between K-injective complexes, therefore so is
$\zeta\set i_*(\sigma\smcirc i^*\xi)$, as is the induced map of cones
$\epsilon\:\Cal C_\eta\to\Cal C_{\zeta\circ\eta\>}$.

{}From the commutative diagram of complexes
$$
\def\1{\F}
\def\2{i_*i^*(\F\>)}
\def\4{i_*\Cal L}
\def\5{\Cal K_\flat\<\<\otimes\F\>}
\def\6{i_*i^*(\Cal K_\flat\<\<\otimes\F\>)}
\CD
\1@>\eta>>\2 @>\zeta>> \4 \\
@V\xi VV @VVi_*i^*\xi V @|\\
\5@>>>\6@>>i_*\sigma>\4
\endCD\tag 3.2.1
$$
we deduce a map of cones
$$
\Cal C_\xi @>>> \Cal C_{\zeta\circ\eta\>}\tag 3.2.2
$$
and hence a composed $\D$-map
$$
\delta(\F\>)\: \KK(\bt,\F\>)\cong \Cal C_\xi[-1]\to\Cal C_{\zeta\circ\eta}[-1]
@>\epsilon^{-1}>>\Cal C_\eta[-1]\cong\vg Z(\F\>)\cong \R\vg Z(\F\>),
$$
easily checked to be functorial in~$\F$.

d) To check that $\delta\smcirc\delta'$ is as asserted above,
``factor" the first square in~(3.2.1) as
$$
\CD
\F @>\pi>>\F/\vgpt(\F\>) @>\eta'>> i_*i^*(\F\>) \\
@. @V\xi'VV @VVV \\
@. \Cal K_\flat\<\<\otimes\F @>>> i_*i^*(\Cal K_\flat\<\<\otimes\F\>)\,,
\endCD
$$
derive the commutative diagram
$$
\CD
\Cal C_\pi[-1] @>\text{via }\eta'>> \Cal C_\eta[-1] \\
@V\text{via }\xi'VV @VV\epsilon\>[-1] V\\
\Cal C_\xi[-1] @>>(3.2.2)> \Cal C_{\zeta\circ\eta}[-1]\,,
\endCD
$$
and using a) and b), identify the $\D$-map labeled ``\,via $\xi'\,$''
(resp.~``\,via  $\eta'\,$'') with~$\delta'$
(resp.~the inclusion map~$\vgpt(\F\>)\hookrightarrow \vg Z(\F\>)$).
\smallskip
The next Lemma is a derived-category version of
\cite{Gr, p.\,20, Prop.\,5} and \cite{H, p.\,98, Prop.\,4.3, b)}
(from which it follows easily if the complex $\F$ is bounded-below
or if the functor $\vg Z$ has finite homological dimension).

\proclaim{Lemma (3.2.3)}
If\/ $(X,\O)$ is a scheme,  $\bt$ is a finite sequence in\/~$\Gamma(X,\O),$\
and if\/ $Z\set\text{\rm Supp}(\O/\bt\O),$\ then\/
$\delta(\F\>)\:\KK(\bt,\F\>)\to\R\vg{Z}(\F\>)$
is an isomorphism for all\/~$\F\in\nmb\D\qcd(X)$.

\endproclaim

\proof
The question is local, so we may assume $X$ to be affine,
say $X=\text{Spec}(R)$. Let $i\:U\set (X\setminus Z)\hookrightarrow X$ be the
inclusion, a quasi-compact map (since $U$ is quasi-compact). Let
$\Cal K_\flat$ be as in the definition of~$\delta$, so that $\Cal
K_\flat=i_*i^*\Cal K_\flat$.
Also, the \v Cech resolution $i^*\xi(\O)\:\O_U\to i^*\Cal K_\flat$ (see c)
above) is
{\it $i_*$-acyclic,} i.e., $R^pi_*(i^*\Cal K_\flat^q)=\nmb0$
for all~$p>0$ and $q\ge 0$: indeed,
$i^*\Cal K_\flat^q$ is a direct sum of sheaves of the form $j_*\O_V$, where
$V\subset U$ is an open set of the form  $\text{Spec}(R_{\>t})$
($t$ a product of some members of~$\bt$)
and $j\:V\hookrightarrow U$ is the inclusion map;
and since $V$ is affine, therefore
$$
i_*(j_*\O_V)=(ij)_*\O_V=\R(ij)_*\O_V=\R\>i_*(\R j_*\O_V)=\R\>i_*(j_*\O_V),
$$
whence $i_*(i^*\Cal K_\flat^q)=\R\>i_*(i^*\Cal K_\flat^q)$.
It follows that $\Cal K_\flat= i_*(i^*\Cal K_\flat)\cong \R\>i_*(\O_U).$

Since the bounded complex $\Cal K_\flat$ is {\it flat,}
we conclude that the bottom row of~(3.2.1) is isomorphic
in $\D$ to the canonical composition
$$
\smash{\R\>i_*\O_U\<\<\Otimes\kern-.8pt \F\>} \to
\smash{ \R\>i_*i^*(\R\>i_*\O_U\<\<\Otimes\kern-.8pt \F\>)}
  \iso\< \smash{\R\>i_*(i^*\R\>i_*\O_U\<\<\Otimes\kern-.6pt i^*\F\>)}
  \iso\<\smash{\R\>i_*(\O_U\<\<\Otimes\kern-.6pt i^*\F\>)}
$$
{\it which composition is an isomorphism
for any\/} $\F\in \D\qcd(X)$.
This instance of the ``projection isomorphism''
of~\cite{H, p.\,106} (where the hypotheses are too restrictive) is
shown in \cite{L, Prop.\,3.9.4} to hold in the necessary generality.
It follows that the map~$\Cal C_\xi\to \Cal C_{\zeta\circ\eta}$ in~(3.2.2) is a
$\D$-isomorphism, whence the assertion.
\endproof

{}From the implication $(1)\Rightarrow (2)$ of Lemma~(3.1.1)---whose proof does
not need $X$ to be concentrated---we now obtain:

\proclaim{Corollary (3.2.4)}
If $Z$ is a proregularly embedded subscheme of the scheme~$X$ then for all
$\F\in\nmb\D\qcd(X),$\
the natural map\/
$\R\vgp Z(\F\>)\to\R\vg Z(\F\>)$ is an isomorphism.
\endproclaim

\proclaim{Corollary (3.2.5)}
Let\/ $(X,\O)$ be a scheme and\/ $Z\subset X$ a closed subscheme such that the
inclusion $(X\setminus Z)\hookrightarrow X$ is quasi-compact.

{\rm(i)} The map
$
\botsmash{\psi\:\E\Otimes\R\vg Z(\F\>)\iso \R\vg Z(\E\Otimes\F\>)}
$
of\/ $(2.1)$ is an isomorphism for all
$\E,\F \in\D\qcd(X)$.\vadjust{\kern.6pt}

{\rm(ii)} If\/ $Z^*\subset X$ is a closed subscheme such that
$(X\setminus Z^*)\hookrightarrow X$ is quasi-compact,
then the the natural functorial map\/
$\R\vg{Z\cap Z^*}\E\to\R\vg Z\R\vg{Z^*}\E$
is an isomorphism for all
$\E\in\D\qcd(X)$.\vadjust{\kern1pt}

{\rm(iii)} $\R\vg Z\bigl(\D\qcd(X)\bigr)\subset \D\qcd(X)$.
\endproclaim

\proof
Since $\psi$ is compatible with restriction to open subsets, we may
assume that $X$ is affine, so that $Z=\text{Supp}(\O/\bt\O)$
for some finite sequence~$\bt$ in~$\Gamma(X,\O)$.
We may also assume that $\E$ is K-flat, and  then check that
the following diagram---whose top row is the natural
isomorphism---commutes:
$$
\CD
\E\otimes\KK(\bt,\F\>) @>\Iso>>\KK(\bt,\E\otimes\F\>) \\
@V \text{via }\delta(\F\>) VV @VV\delta(\E\otimes\F\>)V  \\
\E\otimes\R\vg Z(\F\>) @>\psi>>\R \vg Z(\E\otimes\F\>)
 \endCD
$$
Since both $\E$ and $\F$ are in $\D\qcd(X)$, so is \smash{$\E\Otimes\F$}:
express\vadjust{\kern1pt} $\E$ and $\F$ as \smash{$\drlm\!$'s} of
bounded-above truncations to reduce to where
$\E,\F\in \D\qcd{}^{\mkern-17mu\umi}\medspace$, a\vadjust{\kern.5pt} case
treated in
\cite{H, p.\,98, Prop.\,4.3}.
By Lemma~(3.2.3) the maps $\delta(\F\>)$
and $\delta(\smash{\E\Otimes\F})$ are isomorphisms,
and assertion~(i) of the Corollary follows.

Assertion~(iii) follows at once from~(3.2.3), see \cite{H, p.\,98, Prop.\,4.3}.
And then (ii) follows from~(3.2.3), since
$\KK((\bt,\bt\<^*),-)\<=\<\KK(\bt)\jog\otimes\jog \KK(\bt\<^*,-)$.
\endproof

\eightpoint
\remark{Remark} As might be expected, assertion (ii) in~(3.2.5) holds for
{\it all\/} $\E\in\D(X)$. This is because\- $\R\vg Z$ can be
computed via ``K-flabby" resolutions, and because for any injective
K-injective complex ~$\J\<$, $\vg{Z^*\<}(\J)$ is K-flabby
(see e.g., \cite{Sp, p.\,146, Prop.\,6.4 and p.\,142, Prop.\,5.15(b)}, and use
the natural triangle $\vg{Z^*\<}(\J)\to\J\to j_*j^*\<\J$ where
$j\:(X\setminus Z^*)\hookrightarrow X$ is the inclusion).
\endremark
\par\tenpoint

\proclaim{Proposition (3.2.6)}\!
Let\/ $(X,\O) $ be a quasi-compact separated scheme, and\/
$Z\!\>\subset\nmb\!\> X$ a~
closed subscheme such that $X\setminus Z$ is quasi-compact.
The following are equivalent:
\smallskip
$(1)$ $Z$ is proregularly embedded in\/~$X$.

$(2)$ The natural functorial map\/
$\boldkey j\smcirc\R\varGamma\>\qcu_{\mkern-8mu Z}
 \to \R\vg Z\smcirc \boldkey j\
($see~Proposition\/~$(3.1.8))$ is an isomorphism.

$(3)$ The natural functorial maps\/
$\boldkey j\smcirc\R\varGamma\>\qcu_{\mkern-8mu Z}
 \to \R\vgp Z\smcirc \boldkey j \to \R\vg Z\smcirc \boldkey j$
are both isomorphisms.
\endproclaim

\proof
$(1)\Rightarrow(3)$. If $Z$ is proregularly embedded in\/~$X$
then Proposition (3.1.8) says that
$\boldkey j\smcirc\R\varGamma\>\qcu_{\mkern-8mu Z}\to \R\vgp Z\smcirc \boldkey
j$
is an isomorphism; and (3.1.1)(2)
and~(3.2.3) give that the natural map
$\R\vgp Z \boldkey j(\F\>)\to\R\vg Z\boldkey j(\F\>)$ is an isomorphism for all
$\F\in\D\qcd(X)$.

$(3)\Rightarrow(2)$. Trivial.
\footnote{One could also prove $(1)\Rightarrow (2)$ without
invoking $\vgp Z$, by imitating the proof of~(3.1.8).%
}

$(2)\Rightarrow(1)$.
Let $i\:Y\hookrightarrow X$ be the inclusion of an affine open subset,
so that $Y\setminus Z$ is quasi-compact, whence
$Y\cap\nmb Z=\nmb\text{Supp}(\O_Y/\bt\O_Y)$ for some finite
sequence~$\bt$ in $\Gamma(Y, \O_Y)$ $(\O_Y\set\O|_Y)$; and let us show
for any $\Cal L$ injective in~$\A\qcd(Y)$ that the canonical map
$\vg{Y\cap Z}\Cal L\to\nmb\R\vg{Y\cap Z}(\Cal L\>)$ is an isomorphism, i.e.,
by~(3.2.3), that $H^n\KK(\bt,\Cal L)=0$ for all $n>\nmb0$. Then (1) will
follow,
by the ring-theoretic analogue of the implication $(2)'\Rightarrow(1)$
in~Lemma~(3.1.1), cf.~\cite{Gr, p.\,24, Lemme 9}.

There is a quasi-coherent $\OX$-module $\Cal L'$ with $i^*\Cal L'=\Cal L$,
and an $\A\qcd$-injective $\J\supset \Cal L'\<$.
Then $\Cal L\subset i^*\<\J$ is a direct summand, and
so for any $n>0$, $H^n\R\vg{Y\cap Z}(\Cal L)$ is a
direct summand of
$H^n\R\vg{Y\cap Z}(i^*\<\J)\cong i^*H^n\R\vg Z(\J)$,
which vanishes if (2)~holds.
Thus $\vg{Y\cap Z}\Cal L\iso\nmb\R\vg{Y\cap Z}(\Cal L\>)$, as desired.
\endproof

\eightpoint
\proclaim\nofrills{Corollary (3.2.7)\usualspace}
{\rm (cf.~\cite{Gr, p.\,24, Cor.\,10}).}
For a concentrated scheme\/~$X,$\ the following are equivalent:

$(1)$ Every closed subscheme\/~$Z$ with\/ $X\setminus Z$ quasi-compact
is proregularly embedded.

$(2)$ For every open immersion\/ $i\:U\hookrightarrow X$ with\/ $U$
quasi-compact, and every\/ $\A\qcd$-injective~$\J\<,$\
the canonical map\/ $\J\to i_*i^*\<\J$ is surjective.
\endproclaim

\proof
Assuming (1), to prove (2) we may assume that $X$ is affine. Then
by \cite{Gr, p.\,16, Cor.\,2.11} we have an exact sequence
$$
0\to\vg Z(\J)\to \J \to i_*i^*\<\J \to H^1\R\vg Z(\J) \to 0,
$$
and so Proposition (3.2.6) yields the conclusion.

Now assume (2) holds, so that for any $\A\qcd$-injective~$\J\<$,
any open immersion $j\:Y\to X$ with $Y$~affine, and any quasi-compact open
$U\subset Y\<$, the restriction $\Gamma(Y,\J)\to\nmb \Gamma(U,\J)$ is
surjective---in other~words, $j^*\<\J$ is {\it quasi-flabby\/}
\cite {Ke, p.\,640}.  To prove (1) it suffices, as in proving the implication
$(2)\Rightarrow(1)$ in~(3.2.6), to show that for any
$\A\qcd$-injective~$\J$ and $n>0$, $H^n\R\vg Z (\J)=0$;
and since the question is local it will be enough to show the same for any
quasi-flabby $\J\<$. For~$n=\nmb1$ this results from the above exact sequence,
and for $n>1$ it results from the isomorphism
$H^n\R\vg Z(\J)\iso H^{n-1}\R\>i_*(i^*\<\J)$ \cite{Gr, p.\,16, Cor.\,2.11},
whose target vanishes because
$i^*\<\J$~is quasi-flabby, hence $i_*$-acyclic \cite{Ke, p.\,641, Cor.\,5}.
\endproof
\tenpoint

\subheading{4. Local isomorphisms} This section provides the proofs which are
still missing from the discussion in~\S2.
Proposition~(4.1) is a
$\D(X)$-variant of Theorem 2.5 in \cite{GM, p.\,447\kern.6pt},
giving a local isomorphism of the
homology of~$\R\sHomb(\R\vg Z(\OX\<),-)$ (called in~\cite{GM\kern.6pt} the
{\it local homology\/} of $X$ at $Z$) to the left-derived functors
of completion along~$Z$. (At least this is done for quasi-coherent flat
$\OX$-modules, but as indicated after~(2.2),
Lemma (4.3) guarantees that's enough.)
Proposition~(4.2) allows\- us to conclude that
on an arbitrary quasi-compact separated scheme~$X\<$, these
isomorphisms---defined via local Koszul complexes---patch together to a
global inverse\- for the map $\Phi(\F,\OX\<)$ of~(2.2).

\proclaim{Proposition (4.1)}
Let\/ $(X,\OX\<)$ be a scheme, let\/ $\bt=(t_1,t_2,\dots,t_\mu)$ be a
proregular sequence in\/~ $\Gamma(X,\OX\<)$ $($Definition~$(3.0.1)),$\ and
set\/ $Z\set\text{\rm Supp}(\OX/\bt\OX\<)$.
Then for any\/ {\rm quasi-coherent} flat\/
$\OX$-module $\P$ there is a\/ $\D(X)$-isomorphism
$$
\RsH{}(\R\vg{Z}(\OX\<),\jog\P)\iso
  \underset{^{r>0}}\to\inlm \P/\bt\<^r\P.
$$
\endproclaim

\proof
Let $\P\to\J$ be an injective resolution. By (3.2.3),
$$
\RsH{}(\R\vg{Z}(\OX\<),\jog\P)\cong \sHomb(\KK(\bt),\jog\J)
\cong \inlm \sHomb(\K(\bt\<^r),\jog\J);
$$
and there are natural maps
$$
\aligned
\pi_i\:H^i\inlm \sHomb(\K(\bt\<^r),\jog\J)
 &\to \inlm H^i\sHomb(\K(\bt\<^r),\jog\J)\\
& \>\>@!@!\cong\,\inlm H^i\sHomb(\K(\bt\<^r),\jog\P),
\endaligned\tag 4.1.1
$$
the last isomorphism holding because $\K(\bt\<^r)$ is a bounded complex
of free finite-rank $\OX$-modules.

It follows easily from the definition of~$\K(\bt\<^r)$ that
$$
H^0\sHomb(\K(\bt\<^r),\jog\P)\cong \P/\bt\<^r\P\jog;
$$
and for $i\ne0$, the implication $(1)\Rightarrow(3)$ in Lemma~(3.1.1) gives
$$
 \inlm H^i\sHomb(\K(\bt\<^r),\jog\P)=0.
$$
It suffices then that each one of the maps $\pi_i$ be an isomorphism; and for
that it's enough that for each affine open $U\subset X$, the natural
composition
$$
\aligned
H^i\Gamma\bigl(U,\,\inlm\sHomb(\K(\bt\<^r),\jog\J)\bigr)\!\<\<
 &\,\>\>\cong H^i\inlm\Hom^\bullet\bigl(\K(\bt\<^r)|_U,\jog\J|_U\bigr)\\
&@>\alpha>>\,\inlm H^i\Hom^\bullet\bigl(\K(\bt\<^r)|_U,\jog\J|_U\bigr)\\
&@>\beta>>\,\inlm\Gamma\bigl(U,H^i\sHomb(\K(\bt\<^r),\jog\J)\bigr)
\endaligned\tag 4.1.2
$$
be an isomorphism. (As $U$ varies, these composed maps form a {\it presheaf\/}
map whose sheafification is ~$\pi_{i\>}$.)

To see that $\beta$ is an isomorphism we can (for notational simplicity)
replace $U$ by~$X$---assumed then to be affine, say $X=\text{Spec}(R)$, write
$\Gamma\E$ for $\Gamma(X,\E)$, and note that since
$\Gamma\>\P\to \Gamma\J$ is a quasi-isomorphism (because
$\P$ is quasi-coherent),
and since $\Gamma\K(\bt\<^r)$ is a finite-rank free $R$-complex, therefore
$$
\aligned
H^i\Hom^\bullet\bigl(\K(\bt\<^r),\jog\J\bigr)
 &\cong H^i\Hom_R^\bullet
  \bigl(\Gamma\K(\bt\<^r),\jog\Gamma\J\bigr)\\
&\cong H^i\Hom_R^\bullet\bigl(\Gamma\K(\bt\<^r),\jog\Gamma\>\P\bigr)\\
&\cong \Gamma H^i\sHomb\<\bigl(\K(\bt\<^r),\jog\P\bigr)_{\phantom{R}}\\
&\cong \Gamma H^i\sHomb\<\bigl(\K(\bt\<^r),\jog\J\bigr).
\endaligned\tag 4.1.3
$$

It remains to be shown that $\alpha$ is an isomorphism; and for that we can
apply \cite{EGA, p.\,66, (13.2.3)}.
As above we may as well assume $X$ affine and $U=X$.

For surjectivity of~$\alpha$, it is enough, by {\it loc.\,cit.,}
that for each $i$, the inverse system
$$
E_r\set\Hom^i\bigl(\K(\bt\<^r),\jog\J\bigr)
=\prod_{0\le p\le\mu}
 \Hom\bigl(\Cal K^p(\bt\<^r),\jog\J^{p+i}\bigr) \qquad(r>0)
$$
satisfy the Mittag-Leffler condition (ML): {\it for each $r$ there is an
$s>r$ such that the images of all the maps $E_{s+n}\to E_r\ (n\ge0)$
are the same.} But we have
$$
\Hom\bigl(\Cal K^p(\bt\<^r),\jog\J^{p+i}\bigr)\cong
\prod_\sigma\,\J_{r\<,@,\sigma}
$$
where $\sigma$ ranges over all $p$-element subsets of $\{1,2,\dots,\mu\}$,
and\vadjust{\kern.5pt} $\J_{r\<,@,\sigma}\set\nmb\J^{p+i}$ for all $r$ and
$\sigma$; and for $s>r$, the corresponding map
$\prod_\sigma\,\J_{s\<,@,\sigma}\to \prod_\sigma\,\J_{r\<,@,\sigma}$
is the direct product of the maps $\J_{s\<,@,\sigma}\to \J_{r\<,@,\sigma}$
given by multiplication by~$t_\sigma^{s-r}$ where
$t_\sigma\set\prod_{j\in\sigma}t_j$.
Thus we need only show there is an~$N$ such that
$t_\sigma^{N+n}\J_{r\<,@,\sigma}=t_\sigma^{N}\<\J_{r\<,@,\sigma}$ for all $r$,
$\sigma$, and $n\ge 0$. But $X$ being affine we have the equivalence
$(1)\Leftrightarrow(2)$ in
Lemma~(3.1.1), which implies that any permutation of~$\bt$ is proregular.
Taking $r=\nmb1$ and $i=1$ in Definition~(3.0.1), and applying Lemma~(3.1.2)
with
$\I=\nmb\I'=\nmb(0)$, we find then that for each $r$, $\sigma$, and
$j=1,2,\dots,\mu$, there is an $N_j$ such that for all~$n\ge0$,
$\>\>t_j^{N_j+n}\<\J_{r\<,@,\sigma}=t_j^{N_j}\!\J_{r\<,@,\sigma}$. The desired
conclusion follows, with $N=\sup(N_j)$.

For bijectivity of $\alpha$, it is enough, by {\it loc.\,cit.,} that
for each $i$, the inverse system
$$
H^i\Hom^\bullet\bigl(\K(\bt\<^r),\jog\J\bigr)
\cong \Gamma H^i\sHomb\bigl(\K(\bt\<^r),\jog\P\bigr) \qquad(r>0)
$$
(see (4.1.3)) satisfy (ML). For $i=0$, this is just the system
$\Gamma(\P)/\bt\<^r\Gamma(\P)$, with all maps surjective;
and for $i\ne 0$, the system is, by Lemma~(3.1.1)(3), essentially null.
\endproof

\proclaim{Proposition (4.2)} With\/ $X,$\ $\bt,$\ $Z$ and\/ $\P$ as
in\/ Proposition~$(4.1),$\ let
$$
\Psi=\Psi(\P)\:\RsH{}(\R\vg{Z}(\OX\<),\jog\P)\iso
  \underset{^{r>0}}\to\inlm \P/\bt\<^r\P=\Lambda_Z(\P)=\L\Lambda_Z(\P)
$$
be the isomorphism constructed in proving that Proposition
$($easily seen to be independent of the injective resolution\/~$\P\to\J$
used there\/$)$
and let
$$
\Phi=\Phi(\P,\OX\<)\:\L\Lambda_Z(\P)\lra \RsH{}(\R\vg{Z}(\OX\<),\jog\P)
$$
be as in ~{\rm (2.2).} Then\/ $\Phi=\Psi^{-1}\<,$\
and so\/ $\Phi$ is an isomorphism.
\endproclaim

\proof
Let $\OX@>\chi>>\Cal R$ be a quasi-isomorphism with $\Cal R$ a bounded-below
injective complex,
and let $\theta\:\P\otimes\Cal R\to \J$ be a quasi-isomorphism
with $\J$ an injective complex vanishing in all negative degrees.
The composition
$$
\P=\P\otimes\OX @>1\otimes\chi>> \P\otimes\Cal R @>\theta>> \J
$$
is then an injective resolution of~$\P\<$, which can be used to define~$\Psi$.
Also, $\Phi$ can be represented by the composition
$$
\Lambda_Z(\P)@>(2.2)>>\sHomb(\vgp{Z}(\Cal R),\jog\P\otimes\Cal R)
@>\text{via }\theta>>\sHomb(\vgp{Z}(\Cal R),\jog\J).
$$
So it suffices to show that the following composition---which is
$H^0(\Psi\smcirc\Phi)$---is the identity map of~$\>\inlm \P/\bt\<^r\P$.
$$
\aligned
\inlm \P/\bt\<^r\P
  & @>(2.2)>>H^0\sHomb(\vgp Z(\Cal R),\jog\P\otimes\Cal R) \\
  & @>\text{via }\theta>> H^0\sHomb(\vgp Z(\Cal R),\jog\J) \\
  &@>(3.1.1)>> H^0\sHomb(\KK(\bt,\Cal R),\jog\J)   \\
  & @>\<\text{via }@!\chi\<>> H^0\sHomb(\KK(\bt),\jog\J)  \\\vspace{1\jot}
  & =\!\!=\!\!=\!\!=\!\!\!=H^0\inlm \sHomb(\K(\bt\<^r),\jog\J)
  @>\Iso>(4.1.1)>\,\inlm \P/\bt\<^r\P.
\endaligned\tag 4.2.1
$$
This composition is the sheaf map associated to the composed presheaf map
obtained by replacing, throughout, $\sHomb(-,-)$ by $\Gamma(U,\sHomb(-,-))$
and $\P/\bt\<^r\P$ by $\Gamma(U,\P/\bt\<^r\P)$, with $U$ an arbitrary
affine open subset of~$X\<$. Thus we may replace $X$ by $U$, i.e.,
assume $X$ affine, say $X=\text{Spec}(R)$, and then in~(4.2.1)
we may replace $\sHomb$ by
$\Homb=\Gamma(X,\sHomb)$ and $\P/\bt\<^r\P$ by $P/\bt\<^rP\>$ where
$P$ is the $R$-module~$\Gamma(X,\P)$. Note that the arrow labeled~(4.1.1)
remains an isomorphism after these replacements are made, for then
it factors as
$$
\multlinegap{1pt}
\multline
H^0\inlm\Homb(\K(\bt\<^r),\jog\J)\underset\alpha\to\iso
 \inlm H^0\Homb(\K(\bt\<^r),\jog\J) \iso \\
\inlm H^0\Hom_R^\bullet(\Gamma\K(\bt\<^r),\jog\Gamma\J)
\cong\inlm H^0\Hom_R^\bullet(\Gamma\K(\bt\<^r),\jog\Gamma\>\P)\cong\inlm
P/\bt\<^rP.
\endmultline
$$
(In the proof of Proposition (4.1), $\alpha$ is shown to be an isomorphism;
and for the rest see~(4.1.3) and the remarks preceding it.)

As in Lemma (3.1.1), set
$$
\H_{@,@,\bt\<^r}\Cal R\set\sHom(\OX/\bt^r\OX,\jog\Cal R)
 =\ker\bigl(\Cal K^0(\bt^r\<,\Cal R)\to
\Cal K^1\<\<(\bt^r\<,\Cal R)\bigr).
$$

There are natural commutative diagrams of complexes $(r>0)$:
$$
\CD
\vgp Z(\Cal R) @= \drlm \H_{@,@,\bt\<^r}\Cal R @<<< \H_{@,@,\bt\<^r}\Cal R\\
@V\xi VV @. @VV\eta V \\
\KK(\bt, \Cal R) @= \drlm \K(\bt\<^r, \Cal R)
  @<<<\K(\bt\<^r, \Cal R)\\
\endCD
$$

It will suffice then to show that the composition $bcde\<\<f$ in the following
commutative diagram---an expansion of~(4.2.1), mutatis mutandis---is the
identity map of~$\>\inlm P/\bt\<^rP$.

$$
\def\1{\inlm P/\bt\<^rP}
\def\2{H^0\Homb(\vgp Z(\Cal R),\jog\J)}
\def\3{H^0\Homb(\KK(\bt,\Cal R),\jog\J)}
\def\4{H^0\Homb(\KK(\bt),\jog\J)}
\def\5{\inlm H^0\Homb(\H_{\bt\<^r}\Cal R,\jog\J)}
\def\6{\inlm H^0\Homb(\K(\bt\<^r, \Cal R),\jog\J)}
\def\7{\inlm H^0\Homb(\K(\bt\<^r),\jog\J)}
\def\8{\inlm H^0\Homb(\K(\bt\<^r),\jog\P)}
\CD
\1    \\
@VVfV   \\
\2 @>g>> \5 \\
@V\simeq Ve\text{ (via }\xi^{-1})V @AA h\text{ (via }\eta)A \\
\3 @>d>\undercircle1> \6  \\
@V\simeq V\text{via }\chi V @V\simeq V c\text{ (via }\chi) V  \\
\4 @>\Iso>\alpha>\7 @>\Iso>b> \1
\endCD
$$

Commutativity of subdiagram \circle1 implies that $d$
is an isomorphism.
The map~$g$~is surjective, by
\cite{EGA, p.\,66, (13.2.3)}, since $\J$ is an
injective complex so the \-inclusions
$\H_{\bt\<^r}@!@!\Cal R@!\hookrightarrow\<\H_{\bt\<^{r+1}}\Cal R$ induce
{\it surjections\/}
$\Homb(\H_{\bt\<^{r+1}}@!@!\Cal R,\jog\J)@!@!\twoheadrightarrow\<
  \Homb(\H_{\bt\<^r}@!@!\Cal R,\jog\J)\ (r>\nmb0)$.
Hence $h$ too is surjective.
Also $h$ is {\it injective,} hence bijective, as follows from the
natural exact sequences---with $s>0$ and $\Cal C_s$ as in Lemma (3.1.1)(4):
$$
H^0\Homb(\Cal C_s,\jog\J)
 \to \<H^0\Homb(\K(\bt\<^s, \Cal R),\jog\J)
   \to \<H^0\Homb(\H_{\bt\<^s}\<\Cal R,\jog\J).
$$

Now, given an element
$\smash{\bold p\in\inlm P/\bt\<^rP}$, represented
by a sequence $(p_1, p_2, p_3,\dots)$ where $p_r\in\nmb P$ is such that
$p_{r+1}-\nmb p_r\in\nmb \bt\<^rP$ for all ~$r$,
the idea is to construct an element
$\phi^{\>\bold p}\in\inlm H^0\Homb(\K(\bt\<^r, \Cal R),\jog\J)$
with\vadjust{\kern3pt}

(1) $bc(\phi^{\>\bold p})=\bold p$, and

(2) $h(\phi^{\>\bold p})=gf(\bold p)$.\vadjust{\kern3pt}

\noindent Then we will have
$$
bcde\<\<f(\bold p)=bch^{-1}gf(\bold p)=bc(\phi^{\>\bold p})=\bold p,
$$
as desired, completing the proof.

Here's how to get $\phi^{\>\bold p}$. Let $\tau\:\K(\bt\<^r)\to\OX$
be the map of complexes which is in degree zero the identity map of~$\OX$.
Let $\frak p_r\:\OX\to \P$ be the map taking the global section 1 to $p_r$.
Write
$$
p_{r+1}-p_r=\sum_{i=1}^\mu t_i^rq_i\qquad(q_i\in P).
$$

Then the map $\Cal K^1(\bt\<^r)\to \P$ given by the
$\mu$-tuple~$(q_1,\dots,q_\mu)$ provides a homotopy from the composition
$$
\K(\bt\<^r) @>\quad\tau\quad>> \OX @>\frak p_{r+1}-\frak p_r>> \P
$$
to the zero map. Consequently, if $\phi_r$ is the composed map of complexes
$$
\phi_r\:\K(\bt\<^r)\otimes\Cal R @>\tau\otimes 1>>\OX\otimes\Cal R
@>\frak p_r\otimes 1>>\P\otimes\Cal R @>\ \theta\ >>\J
$$
then $\phi_r$ is homotopic to the composition
$\K(\bt\<^r)\otimes\Cal R\to\K(\bt\<^{r+1})\otimes\Cal R@>\phi_{r+1}>>\J\<$,
and~so the family $(\phi_r)_{r>0}$ defines an element ~$\phi^{\>\bold p}\in
\inlm H^0\Homb(\K(\bt\<^r)\otimes\Cal R,\jog\J)$. For this~
$\phi^{\>\bold p}$, the above conditions (1) and~(2) are easily checked.
\endproof
\goodbreak
And finally:

\proclaim{Lemma (4.3)}
If\/ $X$ is a quasi-compact scheme and\/ $Z\subset X$ is a
proregularly embedded closed subset  then
the functor\/ $\R\sHomb(\R\vg Z(\OX\<),-)\:\D\qcd(X)\to \D(X)$
is bounded above and below.
\endproclaim

\proof
Since $X$ is quasi-compact, the question is local, so we may assume that $X$
is affine and that $Z=\text{\rm Supp}(\OX/\bt\OX\<)$ for some proregular
sequence $\bt=(t_1,\dots,t_\mu)$ in $\Gamma(X,\OX\<)$. Lemma~(3.2.3) gives a
functorial isomorphism
$$
\R\sHomb(\R\vg Z(\OX\<),-)\iso\R\sHomb(\KK(\bold t),-).
$$
For any complex
$\E\in\D(X)$ such that $H^i(\E)=0$ whenever $i<i_0\>$, there is
a quasi-isomorphic injective complex~$\J$ vanishing
in all degrees below~$i_0\>$,
and then since the complex $\KK(\bold t)$ vanishes in all degrees outside the
interval~$[0,\mu]$,
$$
H^i\R\sHomb(\R\vg Z(\OX\<),\E)\cong H^i\sHomb(\KK(\bold t),\J)=0
\quad\text{for all}\ i<i_0-\mu.
$$
Thus the functor $\R\sHomb(\R\vg Z(\OX\<),-)$ is bounded below.

To establish boundedness above, suppose $\F\in\D\qcd(X)$ is
such that $H^i(\F\>)=0$ for~all $i>i_0\>$,  and let us prove that
$H^i\R\sHomb(\KK(\bold t),\F\>)=0$ for all $i>i_0\>$.

By \cite{BN, p.\,225, Thm.\,5.1}, we may assume that $\F$ is actually
a quasi-coherent complex, which after truncation may further be assumed
to vanish in degrees $>i_0\>$. Let
$$
f_n\:\tau^{{\sssize\ge}-n}\F\to\J_n\qquad(n\ge 0)
$$
be the inverse system of quasi-isomorphisms of
\cite{Sp, p.\,133, Lemma 3.7\kern.5pt}, where $\tau$ is the truncation
functor and $\J_n$ is an injective complex vanishing in degrees $<-n$.
Writing $\Gamma(-)$ for $\Gamma(X,-)$, we have, for any
$m\in\Bbb Z$ and $n>\max(m,0)$,
natural isomorphisms
$$
H^{-m}\Gamma(\F\>)\iso H^{-m}\Gamma(\tau^{{\sssize\ge}-n}\F\>)
  \iso H^{-m}\Gamma(\J_n),
$$
the second isomorphism holding because both $\tau^{{\sssize\ge}-n}\F$
and $\J_n$ are $\Gamma$-acyclic complexes.
Further, as in the proof of~\cite{Sp, p.\,134, Prop.\,(3.13},
with $\J=\nmb\inlm \J_n\>$ we have natural isomorphisms
$$
H^{-m}\Gamma(\J)\iso H^{-m}\Gamma(\J_n).
$$
Hence {\it the natural map\/ $H^{-m}\Gamma(\F\>)\to H^{-m}\Gamma(\J)$
is an isomorphism for every\/~$m$.}

Knowing that, we can argue just as in the proof of Proposition~(4.1) to
conclude that the maps $\pi_i$
in~(4.1.1)---with $\F$ in place of~$\P$---are isomorphisms for all $i>i_0\>$,
whence the conclusion.
\endproof

\subheading{5. Various dualities reincarnated}
Theorem~(0.3) leads to sheafified\- generalizations ((5.1.3),
respectively~(5.2.3))
of the {\it Warwick Duality\/} theorem of Greenlees and the
{\it Affine Duality\/} theorem of Hartshorne.
In (5.3) we see how together with Grothendieck Duality, Affine Duality
gives a {\it Formal Duality\/} theorem of Hartshorne.
A similar argument yields the related duality theorem of \cite{L2, p.\,188},
which combines local and global duality.
In (5.4), using ~(0.3) and an \cite{EGA}~theorem on homology and completion,
we establish a long exact sequence of Ext functors, which
gives in particular the {\it Peskine-Szpiro duality sequence\/}~(0.4.3).

\proclaim{Corollary (5.1.1)} Let\/ $X$ be a quasi-compact separated scheme
and\/ $Z\subset\nmb X$ a proregularly embedded closed subscheme. Let\/
$\F\in\D\qcd(X),$ let\/ $\gamma\:\R\vg Z\F\to\F$ be the natural
map, and let $\nu\:\F\to \R\> Q\>\LL Z\F$ correspond to
the natural map\/  $\lambda\:\F\to\LL Z\F$ $($see\/ $(0.4)\text{\rm(a)})$.
Then\/ $\gamma$ and\/ $\nu$ induce\/ {\rm isomorphisms}
$$
\gather
\LL Z\R\vg Z\F \iso \LL Z\F, \tag"(i)" \\ \vspace{1\jot}
\R\vg Z\F \iso\R\vg Z\R\> Q\>\LL Z\F. \tag"(ii)"
\endgather
$$
\endproclaim

\proof
Recall from (3.2.5) that $\R\vg Z\F\in\D\qcd(X)$.
Theorem (0.3) transforms the map~(i) into  the map
$$
\R\sHomb(\R\vg Z\OX,\>\R\vg Z\F\>) @>{\<\text{via\;}\gamma\>\>}>>
 \R\sHomb(\R\vg Z\OX,\>\F\>)
$$
which is, by ~(0.4.2), an isomorphism.\vadjust{\kern.5\jot}

\eightpoint
We could also proceed without recourse
to Theorem~(0.3), as follows.
We may assume, by~(1.1), that $\F$ is flat and quasi-coherent. The
question is local, so we can replace $\R\vg Z\F$ by a complex
of the form $\KK(\bt,\F\>)$
(see (3.2.3)), and then via ~(3.2)(c), $\gamma\:\R\vg Z\F\to\nmb \F$
becomes the natural map $\Cal C_\xi[-1]\to\F$ where $\Cal C_\xi$
is the cone of the map
$\xi\:\F\to \Cal K_\flat\otimes\F$ of~(3.2)(b). Since
$\Cal K_\flat\otimes\F=\nmb\bt(\Cal K_\flat\otimes\F\>)$, therefore
$\Lambda_\bt(\Cal K_\flat\otimes\F\>)\set
\inlm\,\bigl((\Cal K_\flat\otimes\F\>)/
 \bt^n(\Cal K_\flat\otimes\F\>)\bigr)=0,$
and so $\LL Z(\gamma)$ is an isomorphism.
\par\tenpoint
\smallskip

As for (ii): with $\Hom\set\nmb\Hom_{\D(X)}\>$ and $\E\in\D\qcd(X)$,  the
composition
$$
\Hom(\R\vg Z\E, \>\F\>)@>>{\<\text{via\;}\nu\>\>}>
 \Hom(\R\vg Z\E, \>\R\> Q\>\LL Z\F\>) @>\Iso>(0.4)\text{\rm(a)}>
  \Hom(\R\vg Z\E, \>\>\LL Z\F\>)
$$
is an isomorphism: it is the map obtained by applying the functor~
$H^0\R\Gamma(X,\>-)$ to the isomorphism~$\lambda'$ of~Theorem (0.3)(bis).
(Recall that
$\R\vgp Z\E\cong\R\vg Z\E$, (3.2.4)).
Hence ``\kern.5pt via $\nu\>$" is an isomorphism, and so by~(0.4.2) the map
$$
\Hom(\R\vg Z\E, \>\R\vg Z\F\>)\to  \Hom(\R\vg Z\E, \>\R\vg Z\R\> Q\>\LL Z\F\>)
$$
induced by~$\nu$ is also an isomorphism. Taking $\E=\R\> Q\>\LL Z\F$, we
see then that the~map~(ii) has an inverse, so it is an isomorphism.
\endproof
\smallskip
\eightpoint
\noindent{\bf Remark (5.1.2).}
We just saw that $\lambda'$ an isomorphism implies that so is (5.1.1)(ii).
Conversely, to show that $\lambda'$ is
an isomorphism, one can reduce via (0.4.2) and (5.1.1)(i) to where
$\F=\nmb\R\vg Z(\F\>)$, then use (5.1.1)(ii) to get for each open $U\subset X$
that the maps
$$
\nopagebreak
\Hom_{\D(U)}\bigl(\R\vg{Z\cap U}\E|_U, \jog \F|_U[i]\bigr)\to
\Hom_{\D(U)}\bigl(\R\vg{Z\cap U}\E|_U, \jog\L\Lambda_{Z\cap U} \F|_U[i]\bigr)
\qquad(i\in\text{\eightBbb Z})
$$
induced by $\lambda$ are all isomorphisms, so that $\lambda'$ induces
homology isomorphisms.

With the notation and relations given in Remark~(0.4)(d), we find that
the map (5.1.1)(ii) is an isomorphism iff the corresponding map
$\R\GG{\jog\bt}\<\<F\to\R\GG{\jog\bt}\<\L\Lambda_\bt^{\phantom{.}}\< F$ is an
isomorphism for any complex of
$A$-modules; in other words, iff Corollary~(0.3.1) holds.
\par\tenpoint
\medskip
 The next result extends Greenlees's
``Warwick Duality" \cite{Gl, p.\,66, Thm.\,4.1} (where $\G=\O_U\>$,
so that $\text{Ext}_U^n(\G,\>i^*\R\>Q\>\LL Z\F\>)=
 \nmb\Bbb H^n(X, \R\>i_*i^*\R\>Q\>\LL Z\F\>)$
is the ``\kern.4pt local Tate cohomology" of~$\F\>$). As before, $Q$ is the
quasi-coherator.

\proclaim{Proposition (5.1.3)} Let\/ $X$ be a quasi-compact separated scheme,
let\/ $Z\subset X$ be a~proregularly embedded closed subscheme,
and let\/ $i\:U=(X\setminus Z)\hookrightarrow X$ be the inclusion. Then for\/
$\G\in\D\qcd(U)$ and\/ $\F\in\D\qcd(X)$ there are natural isomorphisms
$$
\qquad\text{\rm Ext}_{\<U}^n(\G,\>i^*\R\>Q\>\LL Z\F\>)\iso
\text{\rm Ext}_{\<X}^{n+1}(\R\>i_*\G,\>\R\vg Z\F\>)\qquad(n\in\Bbb Z).
$$
\endproclaim

\proof
Since $\G=i^*\R\>i_*\G$, there is a natural isomorphism \cite{Sp, p.\,147,
Prop.\,6.7, (1)}
$$
\R\Hom_U^\bullet(\G,\>i^*\R\>Q\>\LL Z\F\>)\iso
\R\Hom_X^\bullet(\R\>i_*\G,\>\R\>i_*i^*\R\>Q\>\LL Z\F\>).
\tag $*$
$$
The canonical triangle
$\ \R\vg Z\R\>i_*\G\to \R\>i_*\G\to \R\>i_*i^*\R\>i_*\G@>\upl>>\ $
(see (0.4.2.1))
implies $\R\vg Z\R\>i_*\G=\nmb0$; and $\R\>i_*\G\in\D\qcd(X)$
(see \cite{L, (3.9.2)} for the unbounded case); hence
$$
\R\Hom_{\<X}^\bullet(\R\>i_*\G,\R\>Q\>\LL Z\F\>)\cong
\R\Hom_{\<X}^\bullet(\R\vg Z\R\>i_*\G,\>\F\>)=0
$$
(see 0.4(a)), and the triangle
$\ \R\vg Z\R\>Q\>\LL Z\F\to\R\>Q\>\LL Z\F\to \R\>i_*i^*\R\>Q\>\LL Z\F@>\upl>>\
$
yields a natural\- isomorphism
$$
\R\Hom_{\<X}^\bullet(\R\>i_*\G,\>\R\>i_*i^*\R\>Q\>\LL Z\F\>)\<\iso\<
\R\Hom_{\<X}^\bullet\<\bigl(\R\>i_*\G,\>\R \vg Z\R\>Q\>\LL Z\F\>[1]\bigr).
\tag $*@,@,*$
$$
By (5.1.1)(ii) there is a natural isomorphism
$$
\R\Hom_{\<X}^\bullet\<\bigl(\R\>i_*\G,\>\R \vg Z\R\>Q\>\LL Z\F\>[1]\bigr)\iso
\R\Hom_{\<X}^\bullet\<\bigl(\R\>i_*\G,\>\R \vg Z\F\>[1]\bigr).
\tag $*\!*\!*$
$$
Compose the isomorphisms $(*)$, $(*@,@,*)$, $(*\!*\!*)$, and
take homology to conclude.
\endproof
\smallskip
\eightpoint
\remark{Remark} The complex
$
\Cal T_Z\F\set
\RssH X(\R\>i_*\O_U[-1],\>\R \vg Z\F\>),
$
whose hyperhomology
$$
\text{\eightBbb T}_{\<\<Z}^{\>\bullet}(X\<,\F\>)\set
\text{\eightBbb H}^{\>\bullet}(X\<,\>\Cal T_Z\F\>)\set
 H^\bullet\R\Gamma(X,\>\Cal T_Z\F\>) \underset{(5.1.3)}\to \cong
 \text{\eightBbb H}^{\>\bullet}(U\<,\>i^*\R\>Q\>\LL Z\F\>)
$$
\vskip-5pt\noindent
is the {\it local Tate cohomology\/} of~$\F$, is the summit of a triangle based
on the canonical map
$\RssH X(\O_X,\jog\R\vg Z\F\>)\to\RssH X(\R\vg Z\O_X,\jog\R\vg Z\F\>)$, a map
isomorphic via (0.3) and~(5.1.1)(i)
to the natural composition $\R\vg Z\F\to\F\to\LL Z\F$. So there is a long exact
sequence
$$
\cdots\to\text{\eightBbb H}_Z^{\>n}(X,\F\>)\to\text{\eightBbb H}^{\>n}(X,\LL
Z\F\>)
 \to\text{\eightBbb T}_{\<\<Z}^{\>n}(X,\F\>) \to\text{\eightBbb
H}_Z^{\>n+1}(X,\F\>) \to\cdots
$$
and thus, as Greenlees points out, local Tate cohomology pastes together the
right-derived functors of~$\vg Z$ and the left-derived functors of~$\Lambda_Z$.
\par
\endremark
\tenpoint
\medskip

{\bf (5.2).} Next, we derive a generalized form of
Affine Duality \cite{H2, p.\,152, Thm.\,4.1}, see
Corollary (5.2.3): ``double dual = completion".

\proclaim{Proposition (5.2.1)} Let\/ $X$ be a scheme and $Z\subset X$ a
closed subscheme. Then for any $\E,\F\in\D(X)$ there is a natural isomorphism
$$
\R\vg Z\RsH{}(\E,\>\F\>)\iso\RsH{}(\E,\>\R\vg Z\F\>).
$$

If in addition\/ $X$ is quasi-compact and separated, $Z$ is proregularly
embedded, $\F\in\nmb\D\qcd(X),$\ and\/ $\RsH{}(\E,\>\F\>)\in\D\qcd(X),$\
then there is a natural isomorphism
$$
\LL Z\RsH{}(\E,\>\F\>)\iso\RsH{}(\E,\>\LL Z\F\>).
$$
\endproclaim

\proof
Let $i\:(X\setminus Z)\hookrightarrow X$ be the inclusion. Since $i^*$ has
an exact left adjoint (extension by zero), therefore $i^*$ preserves
K-injectivity, and consequently there is a natural isomorphism
$i^*\RsH{}(\E,\>\F\>)\iso\RsH{}(i^*\E,\>i^*\F\>)$.
The first assertion results then from the commutative diagram,
whose rows are triangles (see (0.4.2.1)):
$$
\minCDarrowwidth{.3in}
\CD
\R\vg Z\RsH{}(\E,\>\F\>)@>>>\RsH{}(\E,\>\F\>)@>>>
\R\>i_*i^*\RsH{}(\E,\>\F\>)@>\upl>> \\
@. @| @V\simeq V\text{\cite{Sp, p.\,147, 6.7}} V\\
\RsH{}(\E,\>\R\vg Z\F\>)@>>>\RsH{}(\E,\>\F\>)@>>>
\RsH{}(\E,\>\R\>i_*i^*\F\>)@>\upl>>
\endCD
$$
\smallskip\noindent
The second assertion is given by the sequence of natural isomorphisms
$$
\spreadlines{1\jot}
\align
\ \L\Lambda_Z\R\sHomb\bigl(\E,\jog\F\>\bigr)
 &\underset{(0.3)}\to
  \iso \R\sHomb\bigl(\R\vg Z\OX,\jog\R\sHomb(\E,\>\F\>)\bigr)\\
\vspace{-1.1\jot}
 &\iso \smash{\R\sHomb\bigl((\R\vg Z\OX)\Otimes\E,\jog\F\>\bigr)}
  \qquad\text{\cite{Sp,~p.\,147, 6.6}}\\
 &\underset{(3.1.6)}\to\iso
  \smash{\R\sHomb\bigl(\R\vgp Z\E,\jog\F\>\bigr)}
   \underset{(0.3)}\to\iso
    \R\sHomb\bigl(\E,\jog\L\Lambda_Z\F\>\bigr).\qquad\square
\endalign
$$
\enddemo
\smallskip

Suppose further that $X$ is noetherian.
Let $\Cal R\in\D\qcd(X)$ have
{\it finite injective dimension\/} \cite{H, p.\,83, p.\,134}.
Then for any $\F\in\Dc(X)$ the complex
$$
\Cal D(\F\>)\set\RsH{}(\F,\>\Cal R)
$$
is in $\D\qcd(X)$ \cite{H, p.\,91, Lemma 3.2 and p.73, Prop.\,7.3},
whence---by (3.2.5)---so~is the ``$Z$-dual" complex
$$
\Cal D_Z(\F\>)\set\R\vg Z\Cal D(\F\>)
\underset{(5.2.1)}\to\cong\RsH{}(\F,\jog\R\vg Z\Cal R).
$$

For example, if $\Cal R$ is a {\it dualizing complex\/} \cite{H, p.\,258},
if $x\in X$ is a closed point, and $\J(x)$ is the
injective $\OX$-module vanishing except at~$x$, where its stalk
is the injective hull of the residue field of the local ring
$\O_{\<\<X\<@!\<,@,@,x}\>$, then by \cite{H, p.\,285},
$$
\Cal D_{\{x\}}(\F\>)=\sHomb\bigl(\F,\>\J(x)\bigr)[-d(x)]
$$
where $d(x)$ is the integer defined in \cite{H, p.\,282}.
\smallskip
As in the proof of the second assertion in~(5.2.1), there is a natural
isomorphism
$$
\LL Z\Cal D(\F\>)=\L\Lambda_Z\R\sHomb(\F,\jog\Cal R)\iso
\R\sHomb(\R\vg Z\F,\jog\Cal R)=\Cal D\R\vg Z(\F\>);
$$
and so if $\F\in\Dc(X)$, whence $\Cal D(\F\>)\in\Dc(X)$, then
there is a natural isomorphism
$$
\LL Z\Cal D\>\Cal D(\F\>)\iso  \Cal D\>\R\vg Z\Cal D(\F\>)=
\Cal D\>\Cal D_Z(\F\>)\underset{(0.4.2)}\to=\Cal D_Z\Cal D_Z(\F\>).
$$
\goodbreak
\noindent Thus:
\proclaim{Corollary (5.2.2)} Let\/ $X$ be a noetherian separated scheme,
let\/ $Z\subset X$ be closed, and let\/ $\Cal R\in \D_{\text{c}}(X)$ have
finite injective dimension. Then for any \/ $\F\in \D_{\text{c}}(X)$ we have,
with preceding notation, canonical isomorphisms
$$
\spreadlines{2pt}
\align
\Cal D\R\vg Z(\F\>)&\iso \LL Z\Cal D(\F\>),\\
\Cal D_Z\Cal D_Z(\F\>)&\iso\LL Z\Cal D\>\Cal D(\F\>).
\endalign
$$
\endproclaim

\proclaim{Corollary (5.2.3)} Let\/ $X$ be
a noetherian separated scheme having a dualizing complex\/~$\Cal R$.
Let $Z\subset X$ be closed, and let $\kappa\:X_{/Z}\to X$ be the
completion map.
Then for\/ $\F\in\nmb\D_{\text{\rm c}}(X)$, and with\/ $\Cal D_Z$ as above,
the~natural map\/ $\beta\:\F\to\nmb\Cal D_Z\Cal D_Z\>\F\>$
factors via an isomorphism
$$
\kappa_*\kappa^*\F\iso \Cal D_Z\Cal D_Z\>\F.
$$
\endproclaim

\proof
Since $\Cal R$ is a dualizing complex,
therefore $\Cal R\in\Dc(X)$, $\Cal R$ has finite injective dimension,
and the natural map $\F\to\Cal D\>\Cal D\F\>$
is an isomorphism \cite{H, p.\,258}.
One checks then that $\beta\>$ factors naturally as:
$$
\qquad\F\to\kappa^*\kappa_*\F\underset{(0.4.1)}\to\iso
\LL Z\F\iso \LL Z\Cal D\>\Cal D\F
  \underset{(5.2.2)}\to\iso \Cal D_Z\Cal D_Z\>\F.
\qquad\square
$$
\enddemo

{\bf(5.3).} Here are some applications of Theorem~(0.3) involving Grothendieck
Duality (abbreviated GD) and basic relations between homology and completion.

Let $A$ be a noetherian local ring, with maximal ideal $m$, and
let $I$ be an injective hull of the $A$-module~$A/m$. Assume
that $Y\set\text{\rm Spec}\>(\<A)$ has a {\it dualizing complex\/}~$\Cal
R_Y@,$,
which we may assume to be {\it normalized\/} \cite{H, p.\,276};
and let $f\:X\to Y$ be a proper scheme-map, so that
$\Cal R_X\set\nmb f^!\Cal R_Y$ is a dualizing complex on~$X$
\cite{\kern.5pt V, p.\,396, Cor.\,3}. For any $\F\in\Dc(X)$,
set
$$
\F\>'\set\Cal D(\F\>)=\RsH{}(\F,\>\Cal R_X)\in\Dc(X).
$$
Let $Z$ be a closed subset of~$f^{-1}\{m\}$, define
$\Cal D_Z(\F\>)$ as in (5.2) to be $\R\vg Z(\F\>')$, and let
$\kappa\:\widehat X\to X$ be the canonical map to~$X$ from its
formal completion along~$Z$.

Hartshorne's Formal Duality theorem \cite{H3, p.\,48, Prop.\,(5.2)}
is a quite special instance
of the following composed isomorphism, for $\F\in\D_{\text{c}}(X)$:\kern1pt
\footnote{Hartshorne requires $Z$, but not necessarily $X$, to be
proper over~$A$. Assuming $f$ separated and finite-type, we can reduce
that situation to the present one by compactifying~$f\/$ \cite{L\"u} .}
$$
\nopagebreak
\alignat2
\ \R\Gamma(\widehat X,\>\kappa^*\F\>)=\R\Gamma(X,\>\kappa_*\kappa^*\F\>)
 &\iso \R\Gamma(X,\>\Cal D_Z\Cal D_Z\>\F\>)
    &&(5.2.3) \\
 &\iso \R\Gamma(X,\>\Cal D\>\Cal D_Z\>\F\>)\qquad\qquad\qquad(5.2.1),\
    &&(0.4.2) \\
 &\,=@!\!\!=\!\!= \R\Hom_{\<\<X}^\bullet(\R\vg Z\F\>'\<,\>\Cal R_X) \\
 &\iso \R\Hom_{\<Y}^\bullet(\R f_*\R\vg Z\F\>'\<,\jog\Cal R_Y) &&\text{
\,(GD)}\\
 &\iso\R\Hom_{\<Y}^\bullet(\R f_*\R\vg Z\F\>'\<,\jog
  \R\vg{\{m\}\<}\Cal R_Y)\qquad
    &&(0.4.2) \\
 &\iso \Hom_A(\R\GG Z\F\>'\<,\>I\>) &&\hskip-15pt\text{\cite{H, p.\,285}}\\
\endalignat
$$
where $\GG Z(-)\set\Gamma\bigl(X,\>\vg Z(-)\bigr)$. The last isomorphism
follows from (0.4.4) because $\w I\cong\R\vg{\{m\}}\<\Cal R_Y$
\cite{H, p.\,285}, and $\R f_*\R\vg Z\F\>'\cong \w{\R\GG Z\F\>'}.\jog$
\footnote
{Some technical points here need attention, especially when
$\F$ is unbounded. First,
GD~holds for unbounded $\F\<$, see \cite {N}.
Next, since $\R\vg Z\F\>'\in\D\qcd(X)$
(3.2.5), therefore $\R f_*\R\vg Z\F\>'\in\nmb\D\qcd(Y)$ \cite{L, (3.9.2)};
and so by (1.3),
$\R f_*\R\vg Z\F\>'\smcong \bigl(\R\Gamma(Y, \R f_*\R\vg Z\F\>')\bigr)^\sim\<$.
Finally, using \cite{Sp, 6.4 and~6.7\kern.5pt} and the fact that $f_*$ and
$\vg Z$ preserve K-flabbiness (see Remark following (3.2.5) above),
one checks that
$
\R\Gamma(Y\<, \R f_*\R\vg Z\F\>')\smcong \R\Gamma(X\<, \R\vg Z\F\>')
 \smcong\R\GG Z\F\>'\<.
$}

Taking homology, we get isomorphisms
$$
\ H^q(\widehat X,\> \kappa^*\F\>)\iso \Hom_A(\text{Ext}_Z^{-q}(\F,\>\Cal
R_X),\>I\>)
 \qquad(\F\in\Dc(X),\;q\in\Bbb Z).\tag 5.3.1
$$
(The functor $\text{Ext}_Z^\bullet$ is reviewed in \S5.4 below).

For example, if $A$ is Gorenstein and $f$ is a
Cohen-Macaulay map of relative dimension $n$, then $\Cal R_Y\cong \O_Y$,
$\Cal R_X\cong \omega[n]$ for some coherent $\O_X$-module~$\omega$ (the
{\it relative dualizing sheaf\/}), and (5.3.1) becomes
$$
H^q(\widehat X, \>\kappa^*\F\>)\iso
\Hom_A(\text{Ext}_Z^{n-q}(\F,\>\omega),\>I\>).
$$
\smallskip
Assume now that $Z=f^{-1}\{m\}$.
For $\F\in\Dc(X)$ the following Lemma (with $J=m$), and the preceding
composition  yield isomorphisms
$$
\R\Hom_{\<X}^\bullet(\F,\>\Cal R_X)\otimes_A\< \hat A
\cong\R\Gamma(X\<,\> \F\>'\>)\otimes_A\< \hat A
\cong\R\Gamma(\widehat X\<,\> \kappa^*\F\>'\>)
\cong\Hom_A(\R\GG Z\F\>''\<,\>I\>).
$$
Thus (since $\F\>''=\F\>)$ there is a natural isomorphism
$$
\qquad\R\Hom_{\<X}^\bullet(\F,\>\Cal R_X)\otimes_A \hat A\iso
 \Hom_A(\R\GG Z\F\<,\>I\>)
\qquad\bigl(\F\in\Dc(X)\bigr).
$$
 Since $\R\Hom_{\<\<X}^\bullet(\F,\>\Cal R_X)$ has noetherian homology modules
therefore $\R\Gamma_{\!Z}\F$ has artinian homology modules,
and Matlis dualization produces a natural isomorphism
$$
\R\GG Z\F \iso \Hom_A(\R\Hom_{\<X}^\bullet(\F,\>\Cal
R_X),\>I\>)\qquad\quad\bigl(\F\in\Dc(X)\bigr).\tag "(5.3.2)\ \;"
$$
For bounded~$\F$, this isomorphism is \cite{L2, p.\,188, Theorem},
deduced there directly from GD and Local Duality (which is the case
$X=Y\<$, $\>f={}$identity map).

\proclaim{Lemma (5.3.3)}
Let\/ $A$ be a noetherian ring, $J$ an\/ $A$-ideal, $\hat A$ the\/
$J$-completion,
$f\:X\to \text{\rm Spec}\>(\<A)$ a finite-type map,
 $Z\set\nmb f^{-1}\text{\rm Spec}\>(\<A/J),$\
and  $\kappa\:\widehat X=\nmb X_{/Z}\to\nmb X$ the canonical flat map.
\smallskip
{\rm(a)} If $\E\in \D\qcd(X)$ has proper support\/
$($i.e., $\E$~is exact outside a subscheme\/~$Y$ of\/~$X$ which
is proper over\/~$\text{\rm Spec}\>(\<A)),$\
then there is a natural isomorphism
$$
\R\Gamma(X,\>\E)\otimes_A\hat A\iso \R\Gamma(\widehat X\<,\> \kappa^*\E).
$$

{\rm(b)} Let\/ $\E\in\Dc(X),$\ $\F\in\D\qcd{}^{\mkern-17mu\upl}\;(X),$\
and suppose either that\/ $\E\in\Dc{}^{\!\!\!\umi}(X)$
or that $\F$ has finite injective dimension.
Suppose further that\/ $\RsH X(\E,\>\F\>)$ has proper support.
Then there is a natural isomorphism
$$
\nopagebreak
\R\Hom_{\!X}^\bullet(\E,\>\F\>)\otimes_A \hat A \iso
\R\Hom_{\!\widehat X}^\bullet(\kappa^*\E,\>\kappa^*\F\>)
$$
Hence, by\/ $(0.3)_{\text{c}}\>,$\ if moreover\/ $\F\in\Dc{}^{\!\!\!\upl}(X)$
then there is a natural isomorphism
$$
\R\Hom_{\!X}^\bullet(\E,\>\F\>)\otimes_A \hat A \iso
\R\Hom_{\!X}^\bullet(\R\vg Z\E,\> \F\>).
$$
\endproclaim

\proof
(a) For bounded-below $\E$, way-out reasoning \cite{H, p.\,68, Prop.\,7.1}
brings us to where $\E=\G$, a single quasi-coherent
$\O_X$-module supported in~$Y\<$.  Since
$\G$ is the \smash{$\drlm$} of its coherent submodules, and
since homology on the noetherian spaces $X$ and $\widehat X$ commutes
with~{$\drlm$},
as does $\kappa^*\<$, we can conclude via \cite{EGA, p.\,129, 4.1.10}.

There is an integer~$d$ such that
$H^n(X\<,\G)=0$ for all $n>d$ and all such~$\G$;
so the same holds for $H^n(\widehat X\<, \kappa^*\G)$, and hence the method
(deriving
from \cite{Sp}) used to prove \cite{L, Prop.\,(3.9.2)} gets us from the bounded
to the unbounded case.\vadjust{\kern2pt}

(b) By (a), and since $\RsH X(\E,\>\F\>)\in \D\qcd{}^{\mkern-17mu\upl}\;(X)$,
\cite{H, p.\,92, 3.3}, it suffices to show that
the natural map $\kappa^*\RsH X(\E,\>\F\>)\to
\RsH {\widehat X}(\kappa^*\E,\>\kappa^*\F\>)$ is an isomorphism.
The question is local, so we can assume $X$ affine
and, $\kappa$ being flat, we can
use \cite{H, p.\,68, Prop.\,7.1} to reduce to the trivial case $\E=\O_X^n$.
\endproof
\goodbreak
\medskip
{\bf(5.4).} The exact sequence (0.4.3) is a special case of the last
sequence in the following Proposition~(5.4.1) (which also generalizes the last
assertion in~(5.3.3(b)).\vadjust{\kern1\jot}

 When $W$ is a locally closed subset of a ringed
space~$X\<$, and $\E,\F\in\D(X)$, then following
\cite{Gr, Expos\'e VI\kern.5pt} one sets
$$
\Ext_W^n\bigl(\E,\>\F\>\bigr)\set
 H^n\<\bigl(\R\Gamma^{\phantom{.}}_{\<\!W}\RsH X(\E,\>\F\>)\bigr)
=H^n\<\bigl(\R\>(\Gamma^{\phantom{.}}_{\<\!W}
 \sHom\!_{X})(\E,\>\F\>)\bigr)\quad(n\in\Bbb Z)
$$
where $\Gamma^{\phantom{\cdot}}_{\<\!W}(-)\set\Gamma\bigl(X,\>\vg W(-)\bigr)$
is the functor of global sections supported in~$W\<$, and the second
equality is justified by \cite{Sp, p.\,146, 6.1(iii) and~6.4} (which
uses the preparatory results 4.5, 5.6, 5.12, and 5.22). It also holds,
via (5.2.1), that
$$
\Ext_W^n\bigl(\E,\>\F\>\bigr)=
 H^n\<\bigl(\R\Hom_X^\bullet(\E,\>\R\vg W\F\>)\bigr).
$$
With $U\set X\setminus W$ there is a canonical triangle (cf.~(0.4.2.1))
$$
\R\Gamma^{\phantom{.}}_{\<\!W}\RsH X(\E,\>\F\>)@>>>
\R\Gamma^{\phantom{.}}_{\<\!X}\RsH X(\E,\>\F\>)@>>>
\R\Gamma^{\phantom{.}}_{\<\!U}\RsH X(\E,\>\F\>)@>\upl>>
$$
whence a long exact cohomology sequence
$$
\cdots\to\Ext^n_W(\E,\>\F\>)\to \Ext^n_{\<\<X}(\E,\>\F\>)\to
 \Ext^n_U(\E,\>\F\>)\to\Ext^{n+1}_W(\E,\>\F\>)\to\cdots
$$
\smallskip
\proclaim{Proposition (5.4.1)} Let\/ $X$ be a noetherian separated scheme,
let\/ $Z\subset X$ be a closed subscheme,
and let\/ $\kappa\:\widehat X= X_{/Z}\to X$ be
the canonical map.
Let $\E\in\D(X)$ and $\F\in\Dc(X)$.
Let $W\subset X$  be closed,
so that\/ $W\cap\nmb Z$ is closed in\/~$\widehat X$.
Then there are natural isomorphisms
$$
\Ext_{W\cap Z}^n(\kappa^*\E,\jog\kappa^*\F\>)\iso
\Ext_{\<\<X}^n(\R\vgp Z\E,\jog\R\vg {W\cap Z}\F\>)\qquad(n\in\Bbb Z),
$$
and so with\/ $U\set X\setminus W$ and\/ $\widehat U\set U_{/Z\cap U}$
there is a long exact sequence
$$
\cdots\to\Ext_{\<\<X}^n(\R\vgp Z\E,\jog\R\vg {W\cap Z}\F\>)\to
\Ext^n_{\<\<\widehat X}(\kappa^*\E,\jog\kappa^*\F\>)\to
 \Ext^n_{\widehat U}(\kappa^*\E,\jog\kappa^*\F\>)\to\cdots
$$
Hence under the assumptions of Lemma\/~{\rm(5.3.3)(b)} there is an
exact sequence
$$
\cdots\to\Ext_{\<\<X}^n(\R\vgp Z\E,\jog\R\vg {W\cap Z}\F\>)\to
\Ext^n_{\<\< X}(\E,\>\F\>)\otimes _A\hat A\to
 \Ext^n_{\widehat U}(\kappa^*\E,\jog\kappa^*\F\>)\to\cdots
 $$
\endproclaim

\proof
There are natural isomorphisms
$$
\align
\kappa_*\R\vg{W\cap Z}\RsH{\widehat X}(\kappa^*\E,\jog\kappa^*\F\>)
 &\iso\R\vg W\kappa_*\RsH{\widehat X}(\kappa^*\E,\jog\kappa^*\F\>) \\
 &\iso\R\vg W\RsH{X}(\R\vgp Z\E,\>\F\>) \\
 &\iso\RsH{X}(\R\vgp Z\E,\>\R\vg W\F\>) \\
 &\iso \RsH{X}(\R\vgp Z\E,\>\R\vg {W\cap Z}\F\>).
\endalign
$$
The first isomorphism results from the equality
$\kappa_*\vg{W\cap Z}=\vg W\kappa_*\>$, since $\kappa_*$
preserves K-flabbiness \cite{Sp, p.\,142, 5.15(b) and p.\,146, 6.4}.
The second comes from $(0.3)_{\text{c}}\>.$ The third comes from~(5.2.1).
The last comes from (0.4.2) and~(3.2.5)(ii).

To conclude, apply the functor
{}~$\R\Gamma^{\phantom{.}}_{\<\!X}$ and take homology.
\endproof

\bye